\documentclass{amsproc}    

\usepackage[numbers]{natbib} 

\usepackage{eurosym}     
\usepackage{amssymb}
\usepackage{amsfonts}
\usepackage{cmtt}
\usepackage{graphicx}
\usepackage{tabularx}
\usepackage{booktabs}
\usepackage{float}
\usepackage[hmargin=2.3cm]{geometry}

\newtheorem{theorem}{\textbf{Theorem}}
\newtheorem{lemma}{\textbf{Lemma}}
\newtheorem{Definition}{\textbf{Definition}}
\newtheorem{rem}{\textbf{Remark}}

\newenvironment{Proof}{{\sc Proof.}}{\vspace{0.0cm} \begin{flushright}$\square$\end{flushright}}

\begin{document}
\let\WriteBookmarks\relax
\def\floatpagepagefraction{1}
\def\textpagefraction{.001}

\title[Secretary Problem with Two Pools of Candidates]{On a Secretary Problem with Two Equal-Sized Pools of Candidates}                      

\author[D.M. Ramsey et al.]{David M. Ramsey$^{1,*}$}
\author[]{Ma{\l}gorzata Kuchta$^{2}$}
\author[]{Micha{\l} Morayne$^{1}$}
\author[]{Ma{\l}gorzata Sulkowska$^{3}$}

\thanks{$^{1}$Faculty of Information and Communication Technology, Wroc\l aw University of Science and Technology, Wroc\l aw, Poland}
\thanks{$^{2}$Faculty of Mathematics, Wroc\l aw University of Science and Technology, Wroc\l aw, Poland}
\thanks{$^{3}$Department of Fundamentals of Computer Science, Wroc\l aw University of Science and Technology, Wroc\l aw, Poland}
\thanks{$^{*}$Corresponding author (ORCID: 0000-0002-7186-1436). E-mail: \texttt{david.ramsey@pwr.edu.pl}}
\thanks{E-mails: \texttt{david.ramsey@pwr.edu.pl}, \texttt{malgorzata.kuchta@pwr.edu.pl},\\ \hspace*{50pt} \texttt{michal.morayne@pwr.edu.pl}, \texttt{malgorzata.sulkowska@pwr.edu.pl}}

%
%
%


\thanks{The authors acknowledge support from the NCN grant: DEC-2015/17/B/ST6/01868.}







\begin{abstract}
This paper presents the solution of a secretary type problem where the decision maker observes a sequence of offers of two types. There are $N$ objects of both types and these objects appear to the decision maker in random order. Objects can only be compared when they are of the same type. The objects of a given type that have already appeared can be ranked according to a complete order. Objects of different types are not comparable. On the appearance of an object, the decision maker must decide whether to accept that object (and thus stop searching) or irrevocably reject it (and thus continue searching). The goal of the decision maker is to accept an object that is the best of its type. The authors present the general form of the optimal policy and the asymptotical solution. The authors give a strong result describing whether an object of relative rank one should be accepted when the same number of objects of both types have been seen. 
\end{abstract}



\keywords{Secretary problem, optimal stopping, dynamic programming, partial order, discrete time stochastic process.}

\maketitle

\section{Introduction}

A lot of research concerning the optimal stopping of processes with discrete time involves generalizations or extensions of the so called secretary problem. In the classical secretary problem \citep{gil1}
 a decision maker (DM) sequentially observes a set of $N$ objects whose values are comparible according to a linear ordering. The DM is allowed to select only one object and selection can only occur at the moment when that object appears. The goal of the DM is to maximize the probability of accepting the absolute best object (of absolute rank one). The objects appear in random order. From the definition of the problem, at monent $n$ the DM only observes the rank of the current object relative to the first $n$ objects to appear. The sequence of absolute ranks is a random permutation of the integers from 1 to $N$. The sequence of the relative ranks of the first $n$ objects is the corresponding permutation based on comparing the first $n$ objects. Hence, these relative ranks form a random permutation of the integers from 1 to $n$. The DM should only consider objects that have a relative rank of 1 and given that the current object has relative rank one, the probability that it is the best of all the objects (has absolute rank one) is equal to the proportion of objects that have already appeared. Since decisions are made sequentially, it is natural to apply dynamic programming to solve such a problem \citep{lin1}. 

Some authors have considered analogous problems in which the DM has greater room for manoeuvre, for example the DM could return to a previously observed object, which is still available with some probability \citep{tam1}, \citep{yan1}. Some authors have considered similar problems with different objectives, for instance the DM wishes to maximize the expected value of the chosen number, or to choose one of the best $k$ objects \citep{sza1}, \citep{van1}. In other cases, the DM has a richer source of information, for example she can observe the precise values of objects from a known distribution, rather than just comparing values \citep{sam1}. Game theoretic versions of the 
secretary problem have also been considered \citep{neu1}, \citep{sak1}.

Another direction of research considers sequential decision problems in which the DM can only compare objects according to a partial ordering and the aim of the DM is to choose one of the maximal elements. In such problems one might consider the optimal policy for specific types of partial orders (e.g. binary trees \citep{kaz1}, \citep{kaz2}, \citep{mor1}) or efficient universal algorithms for more general partial orders \citep{fre1}, \citep{gar-mor}, \citep{geo1}, \citep{koz1}, \citep{pre1}, \citep{sta1}. Similar problems can be defined by assuming that objects have multivariate attributes 
\citep{gne1}, \citep{gne2}. 

A natural extension of the classical secretary problem is to consider two linearly ordered pools of objects of the same cardinality, equal to $N$.  Thus we deal here with a partial ordering consisting of two independent linear orders. It is assumed that these $2N$ objects appear in random order. The DM can only accept an object at the moment of its appearance. The goal of the DM is to choose either of the two best ranked objects.  This problem has turned out to be much more difficult that the original secretary problem and it seems that its full solution has been an open problem for a long time. A number of similar problems have been solved. According to one such model, the pools are interpreted as copies of each other. For example, when there are two pools of $N$ objects, these objects can be interpreted as $N$ sets of identical twins. In such problems, each object is comparable to all of the other objects and only its twin is seen to be equivalent. If the DM has not previously seen a copy of the current object, then she knows that a copy of that object will appear later. Hence, the DM should only ever accept the second copy of an object to appear. Such an object should also be the best of those to have appeared so far \citep{gar1}. This model has been extended to problems in which there are $k$ copies of each object \citep{bay1}. According to another such model, there are $k$ such pools of objects. First, the objects from pool one all appear in random order, before the objects from successive pools appear \citep{kuc2}, \citep{kuc1}. The model considered in this paper differs from those according to which copies of objects appear, since according to the model presented here an object from one pool cannot be compared with any object from the other pool. It differs from models according to which the pools appear in sequence, since according to the model presented here the probability that the next object to appear is from a given pool equals the proportion of the objects that have not appeared so far which come from that pool. In this paper we give the form of the solution to this problem for finite $N$, a procedure to derive this solution and the asymptotic form of the solution as $N\rightarrow \infty$. A (perhaps surprisingly) strong result is presented defining whether an object of relative rank one should be accepted when the same number of objects of both types have been seen. 

Section 2 presents the model. A procedure for deriving the optimal strategy based on dynamic programming is given in Section 3. This section also gives intuitive arguments regarding the general form of the optimal strategy. Section 4 presents the concept of a one step look ahead rule (OSLAR), see \citep{son1}, in the framework of the model presented. The general form of the optimal solution, a threshold rule, is derived in Section 5, which includes a (perhaps surprisingly) strong result regarding the threshold rule for all $N\geq 3$. The asymptotic form of the solution is presented in Section 6. Conclusions and directions for future research are presented in Section 7. A number of the proofs and lemmas are technical. For this reason they are presented in an appendix. The appendix also contains a table describing the 
notation used throghout the paper.

\section{A Mathematical Description of the Model} 

A decision maker (DM) observes a sequence of objects of two types, such that there are $N$ objects of both types. These objects appear in random order and the $n$-th object appears at moment $n$, $n=1,2,\ldots ,2N$. The current object can be ranked in comparison to all of the previous objects of the same type, but cannot be compared to any of the objects of the other type.

The DM can only accept one object and an object can only be accepted at the moment of its appearance. Suppose that the DM is still searching at moment $n$, she can either reject that object (and continue observing the sequence) or accept that object (and thus stop searching). 
The goal of the DM is 
to accept one of the two objects that is the absolute best of its type.  Hence, we define
the reward of the DM to be 1 when such an object is chosen and 0 otherwise. 
The absolute rank of an object is defined to be its rank in comparison to all the objects of its type (the best has absolute rank 1, the worst has absolute rank $N$). The relative rank of an object is defined to be its rank in comparison to the objects of its type that have already appeared.
It follows that the DM should only consider objects that have relative rank 1. Such objects will be referred to as candidates.

Suppose that $i$ objects of type 1 and $j$ objects of type 2 have appeared by moment $n$, where $i+j=n$. Let $\alpha_{i,j}$ denote the probability that the next object is of type 1. From the assumption that the objects appear in random order, it follows that
$\alpha_{i,j} = \frac{N-i}{2N-i-j}$, i.e. the number of remaining objects of type 1 
divided by the total number of objects remaining. Given that the next object is of type 1 (2), then the probability that it has relative rank 1 is $\frac{1}{i+1}$ ($\frac{1}{j+1}$, respectively).
Define a sequence of $2N+1$ four-tuples $T_n =(U_n ,V_n, K_n, R_n)$, 
such that
\begin{description} 
\item[i] $T_0 = (0,0,0,0)$,
\item[ii] $P[(U_{n+1},V_{n+1},K_{n+1})=(U_n +1, V_n,1)]=\alpha_{U_n, V_n}$, 
\item[iii]
$P[(U_{n+1},V_{n+1},K_{n+1})=(U_n , V_n +1,2)]=1-\alpha_{U_n, V_n}$.
\item[iv] $P(R_n = r|U_n = i, K_n=1)=\frac{1}{i}, r=1,2,\ldots i$.
\item[v] $P(R_n = r|V_n = j, K_n=2)=\frac{1}{j}, r=1,2,\ldots j$.
\end{description}
Hence, $K_n$ denotes the type of the $n$-th object to appear, while $U_n$ and $V_n$ denote the number of type 1 and type 2 objects, respectively, to have appeared among the first $n$ objects. Analogously, $R_n$ denotes the relative rank of this object with respect to the objects of the same type that have already appeared. By definition
$U_{2N} = V_{2N}= N$. A realization of the sequence of 4-tuples $(U_n ,V_n,K_n,R_n), 0\leq n\leq 2N$, will be referred to as a trajectory of the search process. 

Given $U_n=i$ and $V_n = j$, the relative ranks of future objects are independent of the
remaining information from the trajectory. Hence, it suffices to define the current state of the process at moment $n$ to be $S(n)=(i,j,k)$ where $i=i(n)$ and $j=j(n)$ are the numbers of objects that have appeared of type $1$ and $2$, respectively.
Define $k$ such that
$k=k(n)=0$ if the current object is not the best of its type seen so far, $k=k(n)=w$ if the current object is the best seen so far (a candidate) of type $w$, $w=1,2$. We may assume that when the DM rejects a candidate in state $(i,j,k)$, $k=1,2$, then the state immediately switches to $(i,j,0)$ before the next object appears.

Denote by $\tau_m$ the moment at which the $m$-th candidate appears, i.e. this is the 
$m$-th smallest value of $n$ for which $k(n)\neq 0$. Let $X_m=1$ when 
the object observed at time $\tau_m$ is the best of all objects of its type, otherwise $X_m=0$. Let $Y_m=E(X_m|\mathcal{F}_{\tau_m})$, where $\mathcal{F}_{\tau_m}$ represents the information given by the states observed at moments $\tau_1 ,\tau_2,\ldots ,\tau_m$. The sequence
$Y_m$ is adapted to the filtration $(\mathcal{F}_{\tau_m})_m$ and is not affected by additional information, e.g. the relative rank of other objects that have appeared so far. Note that when $S(\tau_m)=(i,j,k)$, then $Y_m=i/N$ if $k=1$ and $Y_m=j/N$ if $k=2$.
If $\tau_M$ is the time of the appearance of the final candidate, we formally define that $k(\tau_{M+1})=k(\tau_{M+2})=\ldots=0$ and $Y_{M+1}=Y_{M+2}=\ldots =0$. Thus $(Y_m)_m$
may be considered as a process with infinite time horizon and absorbing state $0$. 
Based on this formulation, the DM obtains a reward of 1 by stopping at either the last decision point at which $K_{\tau_m}=1$ or 
the last decision point at which $K_{\tau_m}=2$.

\section{Derivation of the Optimal Strategy by Dynamic Programming} 
\label{Sec3}

Dynamic programming is a well-known method of solving sequential decision problems when the total reward can be written as a linear combination 
of intermediate payoffs, see \citep{bel1}. This approach can be used to derive the optimal reward and strategy for any particular realization of the problem considered. It
will also be useful in the derivation of general results regarding the general form of the optimal strategy and its asymptotic form. Intuitively, given that $j$ objects of type 2 have been seen, then a candidate of type 1 should be accepted if and only if the number of 
type 1 objects to have appeared is equal to or exceeds some threshold $i_N^{\ast}(j)$, i.e. 
a candidate should be accepted in state $(i,j,1)$ if and only if $i\geq i_N^{\ast}(j)$. It is expected that this threshold is non-increasing in $j$. This result will be proved 
in Section \ref{sec5}.  By symmetry, a  candidate should be accepted in state $(i,j,2)$ if and only if $j\geq i_N^{\ast}(i)$.

Let $v_{i,j}$ denote 
the optimal probability of choosing one of the maximal elements when $i$ objects of type 1 and $j$ objects of type 2 have already been rejected, $0\leq i,j\leq N$. 
This will be called the value function. From the symmetry between the two types, 
$v_{i,j}=v_{j,i}$.
By definition $v_{N,N}=0$. 

Given that the present state is $(i,j,0)$, the next observation is of type 1 with probability $\alpha_{i,j}=\frac{N-i}{2N-i-j}$. In this case, it will be a candidate with probability $\frac{1}{i+1}$ and should be accepted if and only if the probability of such a candidate being the best of all objects of type 1, $\frac{i+1}{N}$, is greater than 
the expected reward from future search $v_{i+1,j}$. Analogously, the next observation is  of type 2 with probability $1-\alpha_{i,j}$. In this case, it will be a candidate with probability $\frac{1}{j+1}$ and should be accepted if and only if the probability of such a candidate being the best of all objects of type 2, $\frac{j+1}{N}$, is greater than the expected reward from future search $v_{i,j+1}$. If an object does not have relative rank 1, then the DM continues searching. Hence, the recursion equation for the value function is given by
\begin{equation}
v_{i,j}  =  \alpha_{i,j} \left[ \frac{iv_{i+1,j}}{i+1} + \max \left\{ \frac{1}{N}, \frac{v_{i+1,j}}{i+1} \right\} \right]
\! + \! (1 \! - \! \alpha_{i,j}) \left[ \frac{jv_{i,j+1}}{j+1} + \max \left\{ \frac{1}{N}, \frac{v_{i,j+1}}{j+1} \right\} \right]. 
\label{recurs}
\end{equation} 

Using this equation, the optimal strategy for a given $N$ can be derived recursively. When the present state is $(i,j,1)$, the current object should be accepted if and only if $\frac{i}{N}\geq v_{i,j}$. 
Similarly, when the present state is $(i,j,2)$, the current object should be accepted if and only if $\frac{j}{N}\geq v_{i,j}$. Note that, by symmetry, if the current object should be accepted in state $(i,j,1)$, then the current object should be accepted in state $(j,i,2)$. Also, when $i\geq j$, if a candidate should be accepted in state $(i,j,2)$, then 
it should be accepted in state $(i,j,1)$.

\subsection{Example} 

Consider a simple example where $N=4$. We have $v_{4,4}=0$, thus the last object to be seen should be accepted. From the state $(i,4,0)$, $i\in \{ 0,1,2,3 \}$, the next object to appear must be of type 1 and applying Equation (\ref{recurs}), it follows that
\begin{equation}
v_{i,4} = \frac{i}{i+1} v_{i+1,4} + \max \left\{ \frac{1}{4},\frac{v_{i+1,4}}{i+1} \right\} . 
\end{equation}
This is analogous to the recursion equation for the classical secretary problem with $N=4$. It follows that $v_{3,4} = \frac{1}{4}, v_{2,4} = \frac{5}{12}$ and
$v_{1,4} = v_{0,4} = \frac{11}{24}$. Comparing $v_{i,j}$ and $\frac{i}{N}$, it can be seen that when all the objects of type 2 have been rejected, a candidate of type 1 should be
accepted if and only if $i\geq i_4^{\ast}(4)= 2$. In general, the threshold 
$i_N^{\ast}(N)$ will be equal to the optimal threshold in the classical secretary problem.

 Comparing $v_{i,j}$ and $\frac{j}{N}$, we obtain the intuitively clear fact that when the last object of type 2 to appear is a candidate, then it should always be acceptable (it is the best of all the type 2 objects with probability 1). 

In a similar way, we can derive the value function and optimal decisions when $j=3,2,1,0$ and $i\leq j$ using the recursion equation
\begin{equation}
v_{i,j} = \! \frac{4-i}{8 \! - \! i \! -\! j} \! \left[ \!\frac{iv_{i+1,j}}{i+1} \! + \! \max \left\{ \frac{1}{4}, \frac{v_{i+1,j}}{i+1} \right\} \right]
\! + \! \frac{4-j}{8 \! -\! i \! - \! j} \! \left[ \! \frac{jv_{i,j+1}}{j+1} \! + \! \max \left\{ \frac{1}{4}, \frac{v_{i,j+1}}{j+1} \right\} \right]. 
\end{equation}
The value function and the optimal strategy for this problem are described in Table 1. The optimal strategy for $N=4$ can be described by the set of thresholds $i_4^{\ast}(j)=3$ for $j=0,1,2,3$ and $i_4^{\ast}(4)=2$. A candidate should be rejected in 
state $(i,i,1)$ when $i\in \{ 1,2\}$, but accepted when $i\in \{ 3,4\}$. Define
$i^{\dagger}_N$ such that $i^{\dagger}_N = \min \{ i: v_{i,i} \leq \frac{i}{N}\}$.  Hence, 
$i^{\dagger}_4 =3$. Note that this is the unique solution of 
$i_4^{\ast}(i)=i$.

A program was 
written in the R language to derive the optimal value function and strategy for a given $N$. 
By inspection, for $2\leq N\leq 20$, $i^{\dagger}_N = \lceil \frac{N+1}{2}\rceil$, where $\lceil x\rceil$ denotes the smallest integer that is not exceeded by $x$ (see also Table 2).
It seems that
$v_{0,0}$ is decreasing in $N$ and for $n\geq 100$, $v_{0,0}\approx 0.5+\frac{1}{2N}$.
 These results will be further explored in Sections \ref{sec5} and \ref{sec6}. 

Considering the possible trajectories of the search process when $N=4$, as long as the DM 
has not yet accepted a candidate, the process will pass into a region in which the DM should 
never reject the next candidate of the type that has appeared more often. For example, 
starting in state $(2,0,0)$ or $(2,1,0)$, then the DM should never reject the next candidate 
of type 1 to appear. By symmetry, starting in state $(0,2,0)$ or $(1,2,0)$, then the DM should never reject the next candidate of type 2 to appear. Moreover, it will eventually become optimal to always accept the next candidate to appear, regardless of its type. For example, starting in the state $(2,2,0)$, the DM should always accept the next candidate to appear. We will prove that these results hold for any $N$ in Section \ref{sec5}.

We now give a simple result, which will be useful in 
deriving the general form of the optimal solution.

\begin{lemma}
 For $N\geq 3$, it is always optimal to reject the first object of either type.
\end{lemma} 

\begin{Proof} We can derive a lower bound on the optimal reward from the search 
procedure by considering the expected reward under the following strategy: ignore any object of type 2 and accept the first object of type 1 that is better than the initial object of type 1 to be seen. 
The probability that the second object of type 1 to be observed is better than the first is 0.5. Given that the second 
object is better than the first, the probability that it is the best of all the objects of type 1 is $2/N$. 
For the third object of type 1 to be accepted under such a strategy, it must be the first object to appear that is better than the first object of type 1, i.e. the second object of type 1 must be worse than the first and then the third must be better than
the first two. The probability of this is $1/2 \times 1/3 = 1/6$. Given that this is true, the 
probability that the third object of type 1 is the best of all the objects of type 1 is $3/N$. Let $u$ be the expected payoff from using this strategy. Using the law of 
total probability, 
\[
u \geq \frac{1}{2} \times \frac{2}{N} + \frac{1}{6} \times \frac{3}{N}>\frac{1}{N}.
\]
Since the probability of the first object of type 1 being the best of all objects of class 1 is 
$1/N$, it follows from the symmetry between the two types that the first object of a given type should never be accepted for $N\geq 3$. 
\end{Proof}

On the other hand, as argued above, when the last object of a given type is a candidate, then it should 
always be accepted.

\begin{table}[h] 
\caption{Value function and optimal strategy (in brackets) for the problem with $N=4$: 0 denotes that 
a candidate of type 1 should be rejected in state $(i,j,1)$, 1 denotes that such a candidate should
be accepted, n.a. denotes that no such decision is applicable. The appropriate actions when a candidate of type 2 appears can be derived from the symmetry of the problem.} \label{tab1}
\begin{tabular*}{0.9\linewidth}{@{\extracolsep{\fill}} llllll @{}} \toprule 
& $j=0$ & $j=1$ & $j=2$ & $j=3$ & $j=4$ \\ \midrule
$i=0$ & 0.62004 (n.a.) & 0.62004 (n.a.) & 0.62176 (n.a.) & 0.56667 (n.a.) & 0.45833 (n.a.) \\ 
$i=1$ & 0.62004 (0) & 0.61875 (0) & 0.61875 (0) & 0.55990 (0) & 0.45833 (0) \\ 
$i=2$ & 0.62176 (0) & 0.61875 (0) & 0.61574 (0) & 0.54861 (0) & 0.41667 (1) \\ 
$i=3$ & 0.56667 (1) & 0.55990 (1) & 0.54861 (1) & 0.43750 (1) & 0.25000 (1) \\ 
$i=4$ & 0.45833 (1) & 0.45833 (1) & 0.41667 (1) & 0.25000 (1) & 0.00000 (1) \\
\bottomrule
\end{tabular*}
\end{table}

\begin{table}[h]
\caption{Optimal value and $i^{\dagger}_N$ for various values of $N$.}
 \label{tab2}
\begin{tabular*}{0.9\linewidth}{@{\extracolsep{\fill}} llllllllll @{}} \toprule 
$N$ & 2 & 5 & 10 & 20 & 50 & 100 & 200 & 500 & 1000 \\ \midrule
$v_{0,0}$ & 0.7500 & 0.5884 & 0.5463 & 0.5231 & 0.5094 & 0.5047 & 0.5024 & 0.5010 & 0.5005 \\ 
$i^{\dagger}_N$ & 2 & 3 &  6 & 11 & 26 & 51 & 101 & 251 & 501
\\ \bottomrule
\end{tabular*}
\end{table}

\section{The One-Step Look-Ahead Rule}
\label{Sec4}

The concept of a one-step look-ahead rule (OSLAR) is often very useful in deriving the general form of the optimal strategy in a wide range of stopping problems. An OSLAR
stops at a particular decision point $\tau_m$, if the expected reward from stopping in the current state is greater than the expected reward from waiting and stopping at the next decision point, $\tau_{m+1}$ (if there is such a decision point).  
In general, such a rule is not necessarily optimal. However, suppose that whenever the condition for stopping according to an OSLAR is satisfied at time $\tau_m$, then it is 
satisfied at time $\tau_{m+1}$. In such a case, the optimal policy is an OSLAR (see Sonin \citep{son1}). In our problem, such an OSLAR would compare the expected reward from accepting a candidate that has just appeared to the expected reward from accepting the next candidate to appear (if no such candidate appears, then the reward is defined to be 0).

The following example shows that such an approach should be used with care in the problem we are considering. Consider a situation in which a candidate of type 1 appears when a number of objects of type 1 have already appeared, but none of type 2 have appeared. In this case, it is very likely that the next candidate to appear will be the first object of type 2 to appear. Intuitively, the expected reward from accepting the current candidate (say the $m$-th) will be greater than the expected reward from accepting the $m+1$-th. However, in this case it is likely that the expected reward from waiting for the $m+2$-th candidate to 
appear will be greater than the expected reward from accepting the $m+1$-th to appear. Hence, the optimal strategy is not an OSLAR.

Let $w_{i,j}$ be the expected reward obtained when $i$ objects of type 1 and $j$ objects of 
type 2 have been already rejected and the DM accepts the next candidate to appear.
It should be noted that whenever the next candidate should be accepted, then $w_{i,j}=v_{i,j}$. By symmetry, $w_{j,i}=w_{i,j}$. In addition, using dynamic programming, as in Equation (\ref{recurs}), 
we obtain
\begin{eqnarray}
w_{i,j} & = & \alpha_{i,j} \left[ \frac{1}{N} + \frac{iw_{i+1,j}}{i+1}\right] + (1-\alpha_{i,j}) \left[ \frac{1}{N} + \frac{jw_{i,j+1}}{j+1}\right] \nonumber \\ 
& = & \frac{1}{N} +\frac{\alpha_{i,j} iw_{i+1,j}}{(i+1)} + \frac{(1-\alpha_{i,j})jw_{i,j+1}}{(j+1)}.
\label{recurs2}
\end{eqnarray}

Note that by definition $w_{0,0}=\frac{1}{N}$, since the first object to appear is always 
a candidate and is the best of its type with probability $\frac{1}{N}$. Similarly, 
$w_{N,0}=w_{0,N}=\frac{1}{N}$. Define
$i^{\bullet}_N$ such that $i^{\bullet}_N = \min \{ i: w_{i,i} \leq \frac{i}{N}\}$.

\begin{lemma}
 For $N\geq 3$, $w_{1,1}>\frac{1}{N}$, i.e. $i^{\bullet}_N \geq 2$. 
\end{lemma} 

\begin{Proof}
Assume that the DM starts in state $(1,1,0)$, i.e. one object of both types has appeared and the second object to appear has just been 
rejected. The next object to appear will be a candidate with probability 0.5. If it is a candidate, it will be the best of its type with probability $\frac{2}{N}$. Suppose that
the next object to appear is not a candidate. Without loss of generality, we may assume that
this object was of type 1. The fourth object to appear will be type 1 with 
probability $\frac{N-2}{2N-3}$. Given it is of type 1, it will be a candidate with probability 
$1/3$ and given it is a candidate, it will be the best of type 1 with probability $\frac{3}{N}$.  The fourth object to appear will be type 2 with 
probability $\frac{N-1}{2N-3}$. Given it is of type 2, it will be a candidate with probability 
$1/2$ and given it is a candidate, it will be the best of type 2 with probability $\frac{2}{N}$.  Using the law of total probability, it follows that
\[
w_{1,1} \geq \frac{1}{2}\times \frac{2}{N} + \frac{1}{2} \left[ \frac{N-2}{2N-3}\times \frac{1}{3}\times \frac{3}{N} + \frac{N-1}{2N-3}\times \frac{1}{2}\times \frac{2}{N}
\right] > \frac{1}{N}.
\]
\end{Proof}

We can find a direct expression for $w_{i,j}$ by considering the probability that the next object to appear is the absolute best of its type. We may assume that $i+j=\tau_m$, i.e. the current moment is a decision point, as the distribution of the state at which the next decision is made depends only on $i$ and $j$. 
Let $Y_m$ be the expected reward from stopping at the $m$-th decision point. 
Assume that the state at the $m$-th decision point is $(i,j,k)$.
From an extension of the law of total probability for any two events $A_1$, $A_2$ such that $P(A_k )>0$ for $k=1,2$
\begin{eqnarray}
\: \: \: w_{i,j} \! = \! E(Y_{m+1}) \! & = & \! E(Y_{m+1}|A_1)P(A_1) \! + \! E(Y_{m+1}|A_2)P(A_2)  +  \label{cand2} \\
& &  +E(Y_{m+1}|A_1^c\cap A_2^c)P(A_1^c\cap A_2^c) - E(Y_{m+1}|A_1 \cap A_2)P(A_1 \cap A_2) . \nonumber
\end{eqnarray}
Let $A_k$ be the event that the overall best object of type $k$ has already appeared, $k=1,2$ and assume that $i,j\geq 1$.
Hence, $P(A_1)=\frac{i}{N}>0 , P(A_2)=\frac{j}{N}>0$ and $E(Y _{m+1}|A_1\cap A_2)=0.$
Now consider the probability that the first candidate of type 1 to appear after the $i$-th object of type 1 is the best of all the objects of type 1.  Denote this probability by $P(Z_i)$. Suppose that the next candidate of type 1 to appear is the $l$-th object of type 1 and this is the best of all objects of type 1. For this to occur, then none of the objects $i+1, i+2,\ldots ,l-1$ of type 1 were candidates and the $l$-th is the best overall. The probability of the first set of events is simply the probability that the best of the first $l-1$ objects of type 1 appeared amongst the first $i$, i.e. $\frac{i}{l-1}$. The probability that the $l$-th is the best object overall from an ordering is simply $1/N$.
Note that the best object overall from an ordering has to be a candidate. Summing over the possible values of $l$, i.e. $i+1 \leq l\leq N$,
\begin{equation}
P(Z_i) = \frac{1}{N}\left( 1 + \frac{i}{i+1} +\frac{i}{i+2} +\ldots +\frac{i}{N-1}\right) \label{cand1}.
\end{equation}

Given $A_1$, if a candidate appears, it must be of type 2. Arguing as in the derivation of Equation (\ref{cand1}),
\begin{equation}
E(Y_{m+1}|A_1)P(A_1)= \frac{iP(Z_j)}{N} = \frac{1}{N}\left(1+\frac{j}{j+1}+\frac{j}{j+2}+\ldots+\frac{j}{N-1}\right)\frac{i}{N}. \label{sklad1}
\end{equation}
Analogously,
\begin{equation}
E(Y_{m+1}|A_2)P(A_2)= \frac{jP(Z_i)}{N} = \frac{1}{N}\left(1+\frac{i}{i+1}+\frac{i}{i+2}+\ldots+\frac{i}{N-1}\right)\frac{j}{N}. \label{sklad2}
\end{equation}
Let $S_k$ be the absolute rank of the best object seen so far of type $k$, $k=1,2$. We have
\[
E(Y_{m+1}|A_1^c\cap A_2^c)P(A_1^c\cap A_2^c)=\sum_{u=2}^{N-i+1}\sum_{v=2}^{N-j+1}E(Y_{m+1}|S_1=u,S_2=v)P[S_1=u,S_2=v].
\]
Note that 
\[
P[S_1=u]=\frac{i}{N}\frac{{N-u\choose i-1}}{{N-1\choose i-1}}, \hspace{.2in} 1\leq u\leq N-i+1,
\]
since for $S_1$ to be equal to  $u$, then the object of type 1 with absolute rank $u$ must appear amongst the first $i$ objects and the remaining $i-1$ objects of type 1 seen so far
must be chosen from the $N-u$ objects of type 1 of absolute rank greater than $u$.
Analogously,
\[
P[S_2=v]=\frac{j}{N}\frac{{N-v\choose j-1}}{{N-1\choose j-1}}, \hspace{.2in} 1\leq v\leq N-j+1.
\]
The next candidate to appear is chosen at random from the $u+v-2$ possible candidates and two of them are the absolute best from their ordering.
From the independence of the ranks of objects of different types, it follows that
\begin{equation}
E(Y_{m+1}|A_1^c\cap A_2^c)P(A_1^c\cap A_2^c)=\frac{ij}{N^2}\sum_{u=2}^{N-i+1}\sum_{v=2}^{N-j+1}\frac{{N-u\choose i-1}}{{N-1\choose i-1}}\frac{{N-v\choose j-1}}{{N-1\choose j-1}}\frac{2}{u+v-2}. \label{sklad3}
\end{equation}

 It follows from Equations (\ref{cand2})-(\ref{sklad3}) that
\begin{equation}
w_{i,j}= \!  \frac{i}{N^2}\! \sum_{r=0}^{N-j-1}\! \frac{j}{j+r}+\frac{j}{N^2}\!\sum_{r=0}^{N-i-1}\! \frac{i}{i+r}
+\frac{ij}{N^2}\! \sum_{u=2}^{N-i+1}\sum_{v=2}^{N-j+1}\frac{{N-u\choose i-1}}{{N-1\choose i-1}}\frac{{N-v\choose j-1}}{{N-1\choose j-1}}\frac{2}{u\! +\! v\! -\! 2}.
\label{OSLAR}
\end{equation}

For completeness, we also derive a formula for $w_{i,0}$ using an analogous 
approach. Note that $w_{0,0}=\frac{1}{N}$. For $i\geq 1$, we have
\begin{equation}
w_{i,0} = E[Y_{m+1}|A_1] P(A_1) + E[Y_{m+1}|A_1^c]P(A_1^c). \label{wi01}
\end{equation}
Given $A_1$, the next candidate to appear must be the first object of type 2 to appear, i.e.
\begin{equation}
E[Y_{m+1}|A_1] P(A_1) = \frac{1}{N} \times \frac{i}{N} = \frac{i}{N^2}.  \label{wi02}
\end{equation}
Arguing analogously to the derivation of Equation (\ref{sklad3}),
\begin{eqnarray}
 E[Y_{m+1}|A_1^c]P(A_1^c) & = & \sum_{u=2}^{N-i+1} E[Y_{m+1}|S_1 =u] P(S_1 =u) \label{wi03} \\
& = &  \frac{i}{N} \sum_{u=2}^{N-i+1} \frac{{N-u\choose i-1}}{{N-1\choose i-1}} \frac{2}{u+N-1} . \nonumber
\end{eqnarray}
It follows from Equations (\ref{wi02}) and (\ref{wi03}) that
\begin{equation}
w_{i,0} = \frac{i}{N^2} + \frac{i}{N} \sum_{u=2}^{N-i+1} \frac{{N-u\choose i-1}}{{N-1\choose i-1}} \frac{2}{u+N-1}. \label{wi05}
\end{equation}

\begin{lemma}\label{monotonicity} 
Suppose that $w_{i,j}\leq\frac{i}{N}$ and $j\geq 1$. It follows that for $i'\geq i$, $w_{i',j}\leq \frac{i'}{N}$.
\end{lemma}
\begin{Proof}  
Suppose that $w^N_{i,j}\leq\frac{i}{N}$, i.e.
$$\frac{i}{N}\geq \frac{i}{N^2}\sum_{r=0}^{N-j-1}\frac{j}{j+r}+\frac{j}{N^2}\sum_{r=0}^{N-i-1}\frac{i}{i+r}
+\frac{ij}{N^2}\sum_{u=2}^{N-i+1}\sum_{v=2}^{N-j+1}\frac{{N-u\choose i-1}}{{N-1\choose i-1}}\frac{{N-v\choose j-1}}{{N-1\choose j-1}}\frac{2}{u+v-2}.$$
Equivalently,
\begin{equation}\label{inequality-one}1\geq \frac{j}{N}\left[\sum_{r=0}^{N-j-1}\frac{1}{j+r}+\sum_{r=0}^{N-i-1}\frac{1}{i+r}
+\sum_{u=2}^{N-i+1}\sum_{v=2}^{N-j+1}\frac{{N-u\choose i-1}}{{N-1\choose i-1}}\frac{{N-v\choose j-1}}{{N-1\choose j-1}}\frac{2}{u+v-2}\right].
\end{equation}
The first sum on the right hand side is independent of $i$. The second sum is clearly decreasing in $i$.
Note that the number of terms in the third sum is decreasing in $i$.
Fixing $u$, $v$ and $j$, the term in this sum is decreasing in $i$ if and only if 
$\frac{{N-u\choose i-1}}{{N-1\choose i-1}}$ is decreasing in $i$. We have
$$\frac{{N-u\choose i-1}}{{N-1\choose i-1}}=\frac{N-u}{N-1}\times \frac{N-u-1}{N-2}\times \ldots \times \frac{N-u-i+2}{N-i+1}.$$
Since each term in this product is less than 1, the product is decreasing in $i$. Hence, the conclusion follows.
\end{Proof}

\begin{lemma} \label{OSLAR2}
For $i\geq i^{\bullet}_N$, $w_{i,i} \leq i/N$. Also, $w_{i,i+1}< \frac{i}{N}$.
\end{lemma} 
\begin{Proof}
From Equation (\ref{recurs2}), we have
\[
w_{i,i} = \frac{1}{N} + \frac{iw_{i,i+1}}{i+1}. 
\]
Assume that $w_{i,i} \leq i/N$. Rearranging the above equation, we obtain
\[
w_{i,i+1} = \frac{i+1}{i} \left[w_{i,i}-\frac{1}{N}\right] \leq \frac{i^2-1}{iN} < \frac{i}{N}.
\]
It follows from Lemma \ref{monotonicity} that $w_{i+1,i+1} < \frac{i+1}{N}$ and Lemma 
\ref{OSLAR2} then follows by induction.
\end{Proof}


\section{General Form of the Optimal Strategy} \label{sec5}

In this section the form of the optimal strategy will be derived using a set of recursion 
arguments. A number of the lemmas are both technical and used as stepping stones to 
prove the main theorem of the paper, which presents the general form of the 
optimal strategy for any $N$. The proofs of these lemmas are given in the appendix. 
The proofs of theorems are given in this section.

Intuitively, whenever a sufficiently large number of objects of both types have appeared, then the DM should accept a candidate of either type.
Hence, we start by deriving the region in which a candidate of either type should be accepted. The following lemma states that when $n$ objects in total have already appeared and the DM accepts the next candidate to 
appear, then the expected payoff of the DM is greatest when the difference between the number of objects of each type to have appeared is minimized.

\begin{lemma} \label{even}
Let $0\leq a\leq b \leq c \leq d\leq N$, such that $a+d=b+c=n$ and $(a,d)\neq (b,c)$. Then $w_{a,d}\leq w_{b,c}.$
\end{lemma}

The second technical lemma is critical in proving that when it is optimal to accept the 
next candidate to appear, then the function $v_{i,j}$ is decreasing in both $i$ and $j$. 

\begin{lemma} \label{w_dec}
Suppose that $j\leq i<N$. It follows that
\[
\frac{iw_{i+1,j}}{i+1} \geq \frac{jw_{i,j+1}}{j+1}, \text{i.e. } \frac{1}{N} + \frac{iw_{i+1,j}}{i+1} \geq \frac{1}{N} + \frac{jw_{i,j+1}}{j+1}.
\]
\end{lemma}

\begin{rem}
Suppose that the DM accepts the next candidate to appear. Lemma \ref{even} can be interpreted as stating that given that the next object is not a candidate,
then the DM prefers the object to be of the type that has appeared less often. Lemma \ref{w_dec} states that when there is no information about the next object,
then the DM prefers the object to be of the type that has appeared more often.
\end{rem}

\begin{theorem} \label{acc_next}
Suppose that $w_{i,j}\leq \min \left\{ \frac{i}{N}, \frac{j
}{N}\right\}$. For $i'\geq i, j'\geq j$, the function 
$w_{i',j'}$ is decreasing in both of its arguments. 
\end{theorem}

\begin{Proof} Note that from symmetry, we may assume that $i\geq j$. It follows from Lemma \ref{monotonicity} that
\[
\frac{1}{N} \geq \frac{w_{i+1,j}}{i+1}; \hspace{.1in} \frac{1}{N} \geq \frac{w_{i,j+1}}{j+1}.
\]
From Equation (\ref{recurs2}), we obtain
\begin{eqnarray*}
w_{i,i} & = & \frac{1}{N} + \frac{iw_{i+1,i}}{i+1} \\
& \geq & \frac{w_{i+1,i}}{i+1} + \frac{iw_{i+1,i}}{i+1} = w_{i+1,i}=w_{i,i+1}. 
\end{eqnarray*}
Suppose that $i>j$. 
From Equation (\ref{recurs2}), we obtain
\[
w_{i,j} = \frac{1}{N}+\frac{i\alpha_{i,j} w_{i+1,j}}{i+1} + \frac{j(1-\alpha_{i,j}) w_{i,j+1}}{j+1}.
\]
From Lemma \ref{w_dec} and the assumptions of the theorem
\[
w_{i,j} \geq \frac{w_{i,j+1}}{j+1} + \frac{j\alpha_{i,j} w_{i,j+1}}{j+1} + \frac{j(1-\alpha_{i,j}) w_{i,j+1}}{j+1} = w_{i,j+1}.
\]
Note that from Lemma \ref{even}, $w_{i+1,j}\leq w_{i,j+1}$ and thus $w_{i+1,j}\leq w_{i,j}$. 

It follows that $w_{i+1,j}\leq \min \left\{ \frac{i+1}{N}, \frac{j}{N}\right\}$
and  $w_{i,j+1}\leq \min \left\{ \frac{i}{N}, \frac{j+1}{N}\right\}$.
Hence, the theorem follows from induction on the number of objects to appear. 
\end{Proof}

\textbf{Corollary to Theorem \ref{acc_next}}: Suppose that 
$w_{i,j} \leq \min \{ \frac{i}{N},\frac{j}{N}\}$. Then it is optimal to 
accept a candidate
in state $(i,j,k)$, $k=1,2$. It is also optimal to accept any candidate that appears at a 
later moment. It follows that $v_{i,j}=w_{i,j}$. Since 
$i^{\bullet}_N = \min \{ i: w_{i,i} \leq \frac{i}{N}\}$ and 
$i^{\dagger}_N = \min \{ i: v_{i,i} \leq \frac{i}{N}\}$, then 
$i^{\bullet}_N =i^{\dagger}_N$.

The condition in the corollary states that it is better to accept a candidate
in state $(i,j,k)$, $k=1,2$  than to accept the next candidate to appear, regardless of the 
type of the next candidate. Theorem \ref{acc_next} guarantees that when the next candidate appears, it is better to accept such a candidate  than to wait for the following candidate to appear. Hence, the corollary follows from the OSLAR property.
It follows that there exists a set, $B^N_3$, of pairs $(i,j)$ such that the first candidate 
to appear after a state $(i,j,0)$ has been reached should always be accepted. 

Note that if, for example, $w_{i,j}>\frac{j}{N}$, then in state $(i,j,2)$ it is better
to reject such a candidate and accept the next candidate to appear than accept the current candidate. It follows that if the inequality $w_{i,j} \leq \min \{ \frac{i}{N},\frac{j}{N}\}$ is not satisfied, then for $i\geq j$ it cannot be optimal to 
accept a candidate in state $(i,j,2)$. Similarly, for $j\geq i$ it cannot be optimal to 
accept a candidate in state $(i,j,1)$.

Intuitively, there will be a set of pairs $(i,j)$ such that $i>j$ and 
$\frac{j}{N}<v_{i,j}\leq \frac{i}{N}$. Thus a candidate should be accepted in 
state $(i,j,1)$, but rejected in state $(i,j,2)$. By symmetry, there will be an analogous set of pairs $(j,i)$ such that a candidate should be accepted in 
state $(j,i,2)$, but rejected in state $(j,i,1)$. Lemma \ref{1N} is useful in both 
describing the boundary of the region $B_3$ with such pairs and showing that if 
it is optimal to stop in state $(i,j,1)$, then it is optimal to stop in
state $(i+1,j,1)$. The proof of this lemma is given in the appendix.

\begin{lemma} \label{1N} For $1\leq i,j \leq N-1$, $|v_{i+1,j}-v_{i,j+1}| \leq \frac{1}{N}$. 
\end{lemma}

\textbf{Corollary to Lemma \ref{1N}}: If a candidate should be accepted in 
state $(i,j,1)$, then a candidate should be accepted in state $(i+1,j,1)$.

\begin{Proof}
Assume that a candidate should be accepted in state $(i,j,1)$, i.e. $v_{i,j}\leq \frac{i}{N}$. We have
\begin{eqnarray*} 
v_{i,j} & \geq & \alpha_{i,j} v_{i+1,j} +(1-\alpha_{i,j} )v_{i,j+1} \\
& \geq & \alpha_{i,j} v_{i+1,j} + (1-\alpha_{i,j})\left( v_{i+1,j}-\frac{1}{N} \right) \\
& \geq & v_{i+1,j} - \frac{1-\alpha_{i,j}}{N}.
\end{eqnarray*}
Hence,
\[
v_{i+1,j} \leq v_{i,j}+\frac{1-\alpha_{i,j}}{N} \leq \frac{i+1-\alpha_{i,j}}{N} \leq \frac{i+1}{N}.
\]
It follows that a candidate should be accepted in state $(i+1,j,1)$.
\end{Proof}

We now consider the form of the boundary of the set $B^N_3$. Assume $i\geq j$. For a fixed $i$, define
$j^{\circ}(i)$ to be the smallest value of $j$ such that $(i,j^{\circ}(i))\in B^N_3$. Note that from the Corollary to Theorem 1, for all $j\geq j^{\circ}(i)$, $(i,j)\in B^N_3$.
 When $i=N$, the problem reduces to the classical secretary problem (CSP). Starting from the state $(N,j,0)$ the DM should accept
a candidate at moment $j'$ if and only if $j' \geq j^{\ast}(N)$, where $j^{\ast}(N)$ is the 
optimal threshold for the CSP. Note that $j^{\ast}(N)$ is non-decreasing in $N$, $j^{\ast}(3)=2$, $j^{\ast}(N)>N\exp(-1)$ and 
$\lim_{N\rightarrow \infty} \frac{j^{\ast}(N)}{N}=\exp(-1)$.  It follows that $(N,j^{\ast}(N)-1)\in B_3$, but $(N,j^{\ast}(N)-2)\notin B_3$. Hence, 
$j^{\circ}(N)=j^{\ast}(N)-1$.

Suppose that $i\geq i^{\bullet}_N$. From the form of $v_{i,j}$,
$j^{\circ} (i)$ is the smallest value of $j$ to satisfy 
$v_{i,j+1}\leq \frac{j+1}{N}$, i.e. if the next object to appear after state $(i,j,0)$ is a candidate of type two, then it should be accepted.  Note that if $(i,j)\in B_3$, then $(i+1,j)\in B_3$.
It follows that $j^{\circ}(i)$ is non-increasing in $i$. In addition, when $v_{i,j}\leq \frac{j}{N}$
 from 
Lemma \ref{1N} it follows that 
$v_{i-1,j+1}\leq \frac{j+1}{N}$. Hence,
whenever $i\geq j+1$ and $(i,j-1)\in B_3$, then $(i-1,j)\in B_3$.
By symmetry, $(i,i-1)\in B_3$, if and only if $(i-1,i)\in B_3$. 
It follows from these facts that for $i\geq i^{\bullet}_N$, $0\leq j^{\circ}(i-1)-j^{\circ}(i)\leq 1$.
Note that for 
$i^{\bullet}_N \leq i\leq N$, we have
 that  $j^{\circ}(N)\leq j^{\circ}(i)\leq j^{\circ}(i^{\bullet}_N)$. From the definition 
of $i^{\bullet}_N$ and the corollary to Theorem 1, $(i^{\bullet}_N,i^{\bullet}_N-1)\in B_3$. It follows that  $(i^{\bullet}_N-1,i^{\bullet}_N)\in B_3$.
 Hence, $j^{\circ}(i^{\bullet}_N)\leq i^{\bullet}_N-1$. Since 
$v_{i^{\bullet}_N-1, i^{\bullet}_N-1}>\frac{i^{\bullet}_N-1}{N}$, it follows from 
Lemma \ref{1N} that $v_{i^{\bullet}_N, i^{\bullet}_N-2}>\frac{i^{\bullet}_N-2}{N}$.
Hence, $(i^{\bullet}_N,i^{\bullet}_N-3)\notin B_3$. It follows that 
$j^{\circ}(i^{\bullet}_N)\geq i^{\bullet}_N-2$.

The following lemma is useful in proving the hypothesis that if a candidate should be accepted in state $(i,j,1)$, then a candidate should be accepted in state 
$(i,j+1,1)$. The proof of this lemma is given in the appendix.

\begin{lemma} \label{1n} Assume that it is optimal to accept an offer in state 
$(i+1,j,1)$, where $j\leq i<N$. It follows that 
\[
\frac{1}{N} + \frac{i v_{i+1,j}}{i+1} \geq v_{i,j+1}.
\]
\end{lemma}

\begin{rem} Lemma \ref{1n} states that when a candidate of type 1 should be accepted it is better for the searcher that the next object to appear is of type 1  rather than an object of type 2 (given that such an object is rejected).
\end{rem}

\textbf{Corollary to Lemma \ref{1n} and Theorem \ref{acc_next}}: Assume that it is optimal to accept an offer in state 
$(i+1,j,1)$, where $i>j$. It follows that $v_{i,j}\geq v_{i,j+1}$. 
Hence, if it is optimal to stop in state $(i,j,1)$, it is also optimal to stop 
in state $(i,j+1,1)$.

This corollary results from the facts that
\begin{description}
\item[i] whenever the next candidate of either type should be accepted, then $v_{i,j}$ is decreasing in both $i$ and $j$.
\item[ii] whenever a candidate should be accepted in state $(i+1,j,1)$, but rejected 
in state $(i,j+1,2)$, then $v_{i,j}$ is a weighted average of 
$\frac{1}{N}+v_{i+1,j}$ and $v_{i,j+1}$, where 
$\frac{1}{N}+v_{i+1,j}\geq v_{i,j+1}$.
\end{description}

Suppose that it is optimal to stop in state $(i,j,1)$. From the corollaries to Lemmas 7 and 8, together with Theorem 1, it follows that it is optimal to stop in both of the 
states $(i+1,j,1)$ and $(i,j+1,1)$. It follows by induction that it is optimal to stop 
in any future state $(i',j',1)$, where $i\leq i'$ and $j\leq j'$. By symmetry, 
if it is optimal to stop in state $(i,j,2)$, then  it is optimal to stop 
in any future state $(i',j',2)$, where $i\leq i'$ and $j\leq j'$.

In summary, the DM should initially reject any candidate, until the trajectory of the search process enters a region in which a candidate of the most often observed type should be accepted. Eventually, the trajectory enters a region in which the first candidate of either type should be accepted. Before giving the general form of the solution, we address the question of the value of $i^{\bullet}_N$ and whether it is possible that the trajectory passes 
directly from the region in which all objects are rejected to the region in which any candidate is accepted. Intuitively, this might occur when the number of objects of 
both type that have been observed is close to $i^{\bullet}_N$.

Since  $v_{i^{\bullet}_N, i^{\bullet}_N} \leq \frac{i^{\bullet}_N}{N}$, a candidate should be accepted in state $(i^{\bullet}_N, i^{\bullet}_N,k)$, where $k\in \{ 1,2\}$ and, from the form of the optimal strategy, at any future decision point.
Also, $v_{i^{\bullet}_N-1, i^{\bullet}_N-1} > \frac{i^{\bullet}_N-1}{N}$. Hence, a candidate should be rejected in state $(i^{\bullet}_N-1, i^{\bullet}_N-1,k)$, where $k\in \{ 1,2\}$ and, from the form of the optimal strategy, at any earlier decision points. 



It follows from Lemma \ref{1N} that when $i\geq j, i+j=2i^{\bullet}_N$, then 
$v_{i,j}\leq \frac{i}{N}$. Hence, a candidate of type 1 should always 
be accepted in state $(i,j,1)$. Analogously, when $i\geq j, i+j=2(i^{\bullet}_N-1)$, then a candidate of type 2 should always be rejected in 
state $(i,j,2)$.

We now consider the optimal actions in state
$(i^{\bullet}_N, i^{\bullet}_N-1,k)$, where $k\in \{1 ,2\}$. The following three cases are possible:

\begin{description}
\item[i] A candidate should be accepted in states $(i^{\bullet}_N, i^{\bullet}_N-1,1)$ and $(i^{\bullet}_N, i^{\bullet}_N-1,2)$.\\ It follows that 
$v_{i^{\bullet}_N, i^{\bullet}_N-1} \leq \frac{i^{\bullet}_N-1}{N}$.
\item[ii] A candidate should be accepted in state $(i^{\bullet}_N, i^{\bullet}_N-1,1)$, but rejected in state $(i^{\bullet}_N, i^{\bullet}_N-1,2)$.\\ It follows that 
$ \frac{i^{\bullet}_N-1}{N}<v_{i^{\bullet}_N, i^{\bullet}_N-1} \leq \frac{i^{\bullet}_N}{N}$.
\item[iii] A candidate should be rejected in states $(i^{\bullet}_N, i^{\bullet}_N-1,1)$ and $(i^{\bullet}_N, i^{\bullet}_N-1,2)$.\\ It follows that 
$v_{i^{\bullet}_N, i^{\bullet}_N-1} > \frac{i^{\bullet}_N}{N}$.
\end{description}
In the first two cases
\begin{equation}
v_{i^{\bullet}_N-1, i^{\bullet}_N-1} = \frac{1}{N}+\frac{(i^{\bullet}_N-1)v_{i^{\bullet}_N, i^{\bullet}_N-1}}{i^{\bullet}_N}.
\label{i-1i-1}
\end{equation}
Hence, in the first case 
\begin{equation}
\frac{i^{\bullet}_N-1}{N} < v_{i^{\bullet}_N-1, i^{\bullet}_N-1} \leq \frac{1}{N}+\frac{(i^{\bullet}_N-1)^2}{Ni^{\bullet}_N} = \frac{i^{\bullet}_N-1}{N} +\frac{1}{Ni^{\bullet}_N} \leq \frac{i_N^{\bullet}}{N}. \label{17}
\end{equation}
Analogously, in the second case
\[
\frac{i^{\bullet}_N-1}{N}+\frac{1}{Ni^{\bullet}_N} < v_{i^{\bullet}_N-1, i^{\bullet}_N-1} \leq \frac{i^{\bullet}_N}{N} .
\] 
In both of these cases, it follows from Lemma \ref{1N} that given $i\geq j, i+j=2i^{\bullet}_N-1$, then the DM should stop in state $(i,j,1)$. Assume that the second case holds. The DM 
should reject a candidate in state $(i,j,2)$. It follows from the symmetry between 
$i$ and $j$, that for $i+j=2i^{\bullet}_N-1$ the DM should accept a candidate of the type that has been more often observed, but reject a candidate of the type that has been observed less often. It follows that in the second case the trajectory cannot pass from the region in which no candidate is accepted directly into the region in which 
any candidate is accepted.

Now assume that the first case holds. When starting from state $(i^{\bullet}_N,i^{\bullet}_N-2,0)$ or $(i^{\bullet}_N-1,i^{\bullet}_N-1,0)$ the DM should always accept the next candidate to appear. It follows from Lemma 
\ref{even} and Equation (\ref{17}) that 
\[
v_{i^{\bullet}_N,i^{\bullet}_N-2}\leq v_{i^{\bullet}_N-1,i^{\bullet}_N-1}
\leq \frac{i^{\bullet}_N}{N}.
\]
Thus the DM should stop in state $(i^{\bullet}_N,i^{\bullet}_N-2,1)$. 
From Lemma \ref{1N} it thus follows that when $i+j=2(i^{\bullet}_N-1)$ and $i\neq j$, then 
a candidate of the type that has been seen more often should be accepted and a candidate of the type that has been seen less often should be rejected. If
$i=j=i^{\bullet}_N-1$, then a candidate should always be rejected.
Hence, the trajectory of the search process can only pass directly from the region in which no candidate is accepted to the region in which any candidate is accepted by passing through a state of the form
$(i^{\bullet}_N-1,i^{\bullet}_N-1, k)$, $k\in \{ 0,1,2\}$.

In the third case, starting in state $(i^{\bullet}_N,i^{\bullet}_N-1,0)$ it is optimal to accept the next candidate to appear. Hence, from the optimality criterion we obtain 
\begin{equation}
v_{i^{\bullet}_N,i^{\bullet}_N-1} =  \frac{1}{N} +\frac{i^{\bullet}_N\alpha_{i^{\bullet}_N,i^{\bullet}_N-1} v_{i^{\bullet}_N+1,i^{\bullet}_N-1}}
{i^{\bullet}_N+1} +\frac{(i^{\bullet}_N-1)(1-\alpha_{i^{\bullet}_N,i^{\bullet}_N-1}) v_{i^{\bullet}_N,i^{\bullet}_N}} 
{i^{\bullet}_N} > \frac{i^{\bullet}_N}{N}. \label{c3}
\end{equation}
From Lemma \ref{w_dec} and Equation (\ref{c3}) it follows that
\begin{eqnarray*}
\frac{1}{N} + \frac{i^{\bullet}_N v_{i^{\bullet}_N+1,i^{\bullet}_N-1}}
{i^{\bullet}_N+1} & > & \frac{i^{\bullet}_N}{N} \\
v_{i^{\bullet}_N+1,i^{\bullet}_N-1} & > & \frac{(i^{\bullet}_N+1)(i^{\bullet}_N-1)}{Ni^{\bullet}_N} > \frac{i^{\bullet}_N-1}{N}.
\end{eqnarray*}
It follows that the DM should reject a candidate in state $(i^{\bullet}_N+1, i^{\bullet}_N-1,2)$. From Lemma \ref{1N} it thus follows that when $i+j=2i^{\bullet}_N$ and $i\neq j$, then 
a candidate of the type that has been seen more often should be accepted and a candidate of the type that has been seen less often should be rejected. If 
$i=j=i^{\bullet}_N$, then a candidate should always be accepted.
Hence, in the third case, the trajectory of the search process can only pass directly from the region in which no candidate is accepted to the region in which any candidate is accepted by passing through a state of the form
$(i^{\bullet}_N,i^{\bullet}_N, k)$, $k\in \{ 0,1,2\}$.

For the first case to occur, the possible range of values of $v_{i_N^{\bullet}-1,i^{\bullet}_N-1}$ is relatively small compared to the second case, particularly when 
$N$ is large. From Inequality (\ref{c3}), in the third case the possible range of values of 
$v_{i_N^{\bullet},i^{\bullet}_N}$ is relatively small. From the numerical calculations made for $N\geq 3$, it was always optimal to stop in state 
$(i^{\bullet}_N,i^{\bullet}_N-1,1)$ and continue in state 
$(i^{\bullet}_N,i^{\bullet}_N-1,2)$. The realization of the problem when $N=2$ is 
specific: $i^{\bullet}_2 = 2$, it is optimal to stop in state $(2,1,1)$ and the 
DM is indifferent between stopping and continuing in state~$(2,1,2)$. 

\begin{theorem} For all $N$, $i_N^{\bullet}=m(N) =\lceil \frac{N+1}{2}\rceil$.
\label{npol}
\end{theorem}

The proof of this theorem is technical and is based the following two lemmas.

\begin{lemma}
For all $N\geq 3$, $w_{m(N)-1,m(N)-1} > \frac{m(N)-1}{N}$. \label{npol1}
\end{lemma}

\begin{lemma}
For all $N\geq 3$, $w_{m(N),m(N)} \leq \frac{m(N)}{N}$. \label{npol2}
\end{lemma}

From Lemma \ref{npol1}, it follows that $i^{\bullet}_N \geq \lceil \frac{N+1}{2} \rceil$.
Theorem \ref{npol} then follows from Lemma \ref{npol2} and Theorem \ref{acc_next}.
The proofs of these lemmas are in the appendix.


Combining the results obtained above, we can state the main hypothesis of the 
paper describing the general form of the optimal strategy. Before doing this, we define 
an auxiliary trajectory, $Q_n$, of the trajectory  $T_n~=~(U_n ,V_n, M_n, R_n)$.

\begin{Definition}
Let $Q_n \! = \! (U'_n , \! V'_n, \! K'_n, \! R_n)$, where $U'_n \! =\! \max\{ \! U_n ,\! V_n\}, \! V'_n \! =\! \min\{ \! U_n , \! V_n\}$,
$K'_n =1$ if the current object is of the type that has been observed more often 
(if $U_n = V_n$, then define $K_n =1$), otherwise $K'_n=2$. As~before, $R_n$ is
the relative rank of the $n$-th object to be observed compared to the previously observed objects of the same type. Denote the sequence of realizations of $U'_n$ and $V'_n$, $n=1,2,\ldots ,2N$ by 
$\{ (i_n ,j_n )\}_{1\leq n\leq 2N}$, where $(i_1 ,j_1 )=(1,0)$ and $i_{2N}=j_{2N}=N$. 
\end{Definition}

\begin{theorem} \label{postac}
Any auxiliary trajectory can be split into three 
subtrajectories: i) $Q_1 ,\ldots, Q_{\mu -1}$, ii) $Q_{\mu },\ldots,$ $Q_{\lambda -1}$ and iii) $Q_{\lambda}, Q_{\lambda +1}, \ldots ,Q_{2N}$, such
is always optimal to continue when $n<\mu$. For $\mu \leq n<\lambda$, it is optimal to accept a candidate of the type that has appeared most often when such a candidate
appears at moment $n$. For $n\geq \lambda$, it is optimal for the DM to accept 
a candidate of either type. The random time $\mu$ denotes the moment at which it first becomes optimal to accept a candidate of the type that has been observed more often. Similalrly, $\lambda$ denotes the moment at which it first becomes optimal to 
accept a candidate of either type. The random time $\mu$ corresponds to the 
first $n$ that satisfies $v_{i_n,j_n}\leq \frac{i_n}{N}$. The random time 
$\lambda$ corresponds to the 
first $n$ that satisfies $v_{i_n,j_n}\leq \frac{j_n}{N}$. The optimal strategy 
of the DM is to reject any candidate when $n<\mu$, accept only a candidate of the 
type to have appeared more often when $\mu \leq n<\lambda$ and 
accept a~candidate of either type when $n\geq \lambda$. 
\end{theorem}

\begin{rem} 
Note that by definition $\mu \leq \lambda$.
\end{rem}

\textbf{Corollary to Theorem \ref{postac}}:
Suppose that $i$ objects of type 1 and $j$ objects of type 2 have appeared, 
$0\leq i,j\leq N$. The pairs $(i,j)$ can be partitioned into the following four sets:
\begin{description}
\item[i] The set of pairs $(i,j)$ such that any object should be rejected in state $(i,j,k)$, $k\in \{ 0,1,2\}$. This set will be denoted by $C^N_0$. For $(i,j)\in C^N_0$, 
$v_{i,j} > \max \{ \frac{i}{N},\frac{j}{N} \}$.
\item[ii] The set of pairs $(i,j)$ such that an object should be accepted in state $(i,j,1)$, but rejected in states $(i,j,0)$ and $(i,j,2)$. This set will be denoted by $C^N_1$.
Thus only candidates of type 1 are accepted in this region.  For $(i,j)\in C^N_1$,
$i>j$ and 
$\frac{j}{N}<v_{i,j} \leq \frac{i}{N}$.
\item[iii] The set of pairs $(i,j)$ such that an object should be accepted in state $(i,j,2)$, but rejected in states $(i,j,0)$ and $(i,j,1)$. This set will be denoted by $C^N_2$. Thus only candidates of type 2 are accepted in this region. For $(i,j)\in C^N_2$,
$j>i$ and 
$\frac{i}{N}<v_{i,j} \leq \frac{j}{N}$.
\item[iv] The set of pairs $(i,j)$ such that an object should be accepted in states $(i,j,1)$ and $(i,j,2)$, but rejected in state $(i,j,0)$. This set will be denoted by $C^N_3$.
Thus candidates of either type  are accepted in this region.  For $(i,j)\in C^N_3$, 
$v_{i,j} \leq \min \{ \frac{i}{N},\frac{j}{N} \}$.
\end{description}

\begin{rem}
All trajectories start in $C^N_0$ and pass at moment, $\mu$, into either $C^N_1$ (when $i_{\mu} >j_{\mu}$)  or $C^N_2$ (when  
$i_{\mu}<j_{\mu}$). At moment $\lambda$, where $\lambda \geq \mu$, the trajectory passes into $C^N_3$ and remains in this set. If $\lambda = \mu$, then the trajectory passes directly from $C^N_0$ to $C^N_3$.
\end{rem}
 
 \begin{rem}
For convenience, it is assumed that any artificial object observed in the artifical states 
$(0,0,k)$, $(i,0,2)$ or $(0,j,1)$ is rejected.
\end{rem}

\begin{rem}
The set $C^N_3$ is very similar, but not identical, to the set $B^N_3$. If $(i,j)\in B^N_3$, then starting in state $(i,j,0)$ it is always optimal to stop when the next candidate appears. However, there will be cases where it is not optimal to stop in state $(i,j,1)$ or in 
state $(i,j,2)$. If $(i,j)\in C^N_3$, then starting in state $(i,j,0)$, it is always optimal to accept the next candidate to appear. Hence,  
$C^N_3 \subset B^N_3$. Also, when $(i,j)\notin C^N_3$, then $(i,j)\in B^N_3$ if and only if $(i+1,j)\in C^N_3$ and $(i,j+1)\in C^N_3$.
\end{rem}

The form of the optimal strategy is illustrated in Fig. 1.

\section{Asymptotic form of the Optimal Strategy} \label{sec6}

Let $v(s,t)$ be the expected reward when a proportion $s$ of objects of type 1 and a proportion $t$ of objects of type 2 have been rejected. 
Define $(s,t,k)$ 
to be the state of the process, where $k=0$ if the present object is not a candidate and when $k=1,2$ then the present object is a candidate of type $k$. 
We wish to derive an integral equation from which we can derive the value function $v$. 
According to the form of the optimal strategy for finite $N$, a candidate of type 1 should be accepted in state $(i,j,1)$ when $i\geq i^{\ast}(j)$. When $N\rightarrow \infty$, this corresponds to  a candidate of type 1 being accepted in state $(s,t,1)$ when $s\geq s^{\ast}(t)$. Analogously, a candidate of type 1 should be accepted in state $(s,t,2)$ if and only if $t\geq t^{\ast}(s)$, where $s^{\ast}\equiv t^{\ast}$.

 By 
symmetry, we have $v(s,t)=v(t,s)$, thus without loss of generality we can assume that $s\leq t$. First, we derive $v(s,t)$ for the set of states $(s,t)$ such that the decision maker 
should accept the next relatively best object to appear, regardless of which type it is. In the limit, the rate at which items of type 2 appear is 
$(1-t)/(1-s)$ times the rate at which items of type 1 appear, i.e. the ratio of the proportion of type 2 objects that have not appeared to the proportion of type 1 objects that have not appeared. Hence, starting from the state $(s,t,0)$, we can define the time $u$ to be the proportion of objects of type 1 that
have appeared. Suppose that the next candidate (of either type) appears at time $u$. From the definitions made above, the state at this point is given by
\begin{equation}
\left[ u, t+\frac{(u-s)(1-t)}{1-s}, k \right] , \label{time}
\end{equation}
where $k=1$ if this relatively best object is of type 1, otherwise $k=2$.

From Ferguson \citep{fer1}, the density function of the time $u$ at which the next relatively best object of type 1 appears is given by $f_{1} (u;s) = s/u^2$, for $s\leq u\leq 1$. 
Note that the event that no relatively best object of type 1 appears after time $s$ is equivalent to the event that the absolutely best object of type 1 appears before $s$, i.e. 
the probability of this event is $s$. Based on the appropriate scaling of time, the density function of the time at which the next relatively best object of type 2 appears is 
 obtained by multiplying the density function of the time at which a candidate of type 2 appears according to the scaled time, $f_1 \left( t+\frac{(u-s)(1-t)}{1-s};t\right)$, by the scaled rate of appearance of such candidates, i.e.
\[
f_2 (u;t) = f_1 \left(t+\frac{(u-s)(1-t)}{1-s};t\right) \frac{1-t}{1-s}= \frac{t(1-t)(1-s)}{[t(1-s)+(u-s)(1-t)]^2}.
\]

From the above considerations, given that the next relatively best object of type 2 (type 1, respectively) appears at time $x$, the probability that it is the next candidate of either type to appear [i.e. no relatively best object of the other type has appeared in 
the time interval $(s,x)]$ is given by $p_2 (x;s,t)$ ($p_1 (x;s,t)$, respectively), where
\[
p_2 (x;s,t) = \frac{s}{x}; \hspace{.1in} p_1 (x;s,t) = \frac{t}{t+\{(x-s)(1-t)/(1-s)\} } = \frac{t(1-s)}{t(1-s)+(x-s)(1-t)}.
\]
The density function of the time to the appearance of a relatively best object from either ordering is given by $f(x;s,t)$, where 
$f(x;s,t)=f_{1} (x;s,t) p_1 (x;s,t) + f_{2} (x;s,t) p_2 (x;s,t)$. Define $r_1 (x;s,t)=f_{1} (x;s,t) p_1 (x;s,t)$ and $r_2 (x;s,t)=f_{2} (x;s,t) p_2 (x;s,t)$. 

Given that the next relatively best object appears at time $u$, if it is of type 1 (type 2, respectively), it is the absolutely best object of type 1
(of type 2, respectively) with probability $u$ [$t+\frac{(u-s)(1-t)}{1-s}$, respectively]. Given that the appearance time of the next relatively best object is $u$, the probability that it 
 is of type 1 is $\frac{r_1 (u;s,t)}{r_1 (u;s,t)+r_2 (u;s,t)}$. It follows that on the region where the next relatively best object should be accepted, the value function 
is given by
\begin{eqnarray}
v(s,t) & \! \! = \! \! & \! \int_s^1 \! \! \left[ \frac{ur_1 (u;s,t)}{r_1 (u;s,t)\! +\! r_2 (u;s,t)}\! + \! \frac{[t(1\! -\! s)\! +\! (u\! -\! s)(1\! -\! t)]r_2 (u;s,t)}{(1-s)[r_1 (u;s,t)\! +\! r_2 (u;s,t)]} \right] \! f(u;s,t) du \nonumber \\
& \! \! = \! \! & \int_s^1 \left[ ur_1 (u;s,t) + \frac{t(1-s)+(u-s)(1-t)}{1-s} r_2 (u;s,t) \right] \nonumber du \\
& \! \! = \! \! & \int_s^1 \left[ \frac{st(2-s-t)}{u[t(1-s)+(u-s)(1-t)]} \right] du \label{int1}.
\end{eqnarray}

First, we consider the case $s=t$. From Equation (\ref{int1}), we obtain
\[
v(s,s) = \int_s^1 \frac{2s^2}{u^2} du = 2s(1-s). 
\]
Suppose a relatively best object of either type appears in state $(s,s,k)$, $k=1,2$. Since such a candidate is the absolute best of its type with probability $s$, it should be 
accepted if and only if $s \geq 2s(1-s)$, i.e. $s\geq 1/2$.

When $s<t$, decomposing the integrand in Equation (\ref{int1}) into partial fractions and integrating, we obtain
\[
v(s,t) = \frac{st(2-s-t)[\ln (t) -\ln (s)]}{t-s}. 
\]
If a relatively best object should be accepted in state $(s,t,1)$, then it should be accepted in state $(s,t,2)$. It follows that any relatively best object should be 
accepted if and only if 
\[
s \geq  \frac{st(2-s-t)[\ln (t) -\ln (s)]}{t-s}. 
\]
It follows that for $s<t$, the boundary of the region in which any relatively best object should be accepted is given by the implicit function
\begin{equation}
1 = \frac{t(2-s-t)[\ln (t) -\ln (s)]}{t-s}. \label{bound1}
\end{equation}
Hence, for a fixed $t$, $t>0.5$, then a candidate should be accepted in state 
$(s,t,1)$, when $s\geq s^{\ast}(t)$, where $s^{\ast} (t)$ satisfies Equation 
(\ref{bound1}). Since $s<t$, a candidate should also be accepted in state $(s,t,2)$.

Now suppose that $(s,t)$ is in the region where only relatively best objects of type 2 should be accepted. In~order to obtain an integral equation for $v(s,t)$, we consider the 
following two events: 

\begin{description}
\item[1.] No relatively best object of type 2 appears before relatively best objects of type 1 become acceptable.
\item[2.] A relatively best object of type 2 appears before relatively best objects of type 1 become acceptable.
\end{description} 

Starting from the state $(s,t,0)$, given that no relatively best object of type 2 appears, in the limit the process will enter the region in which any relatively best object will be 
accepted in the state $[k(s,t),t+\frac{[k(s,t)-s](1-t)}{1-s},0]$, where $v\left( k(s,t),t+\frac{[k(s,t)-s][1-t]}{1-s}\right) =k(s,t)$. The probability that no such relatively best object appears is $\frac{t(1-s)}{t(1-s)+[k(s,t)-s](1-t)}$. 
Now we scale time such that the proportion of objects of type 2 that have appeared at time $t$ is simply $t$. Considering the density function of the time $x$ at which the next 
relatively best object of type 2 appears, since the probability that such an object is the absolute best of the objects  of type 2 is $x$, we obtain
\begin{eqnarray}
v(s,t) & = & \int_t^{\frac{t(1-s)+[k(s,t)-s](1-t)}{1-s}} \frac{t}{x} dx + \frac{t(1-s)k(s,t)}{t(1-s)+[k(s,t)-s](1-t)} \nonumber \\
& = & t\ln \left[ \frac{t(1-s)+[k(s,t)-s](1-t)}{t(1-s)} \right] + \frac{t(1-s)k(s,t)}{t(1-s)+[k(s,t)-s](1-t)}. \label{reg2}
\end{eqnarray}
The lower boundary of the region in which only objects of type 2 are accepted satisfies the equation $t=v(s,t)$. Hence, this lower boundary is given by the implicit function
\begin{equation}
1 = \ln \left[ \frac{t(1-s)+[k(s,t)-s](1-t)}{t(1-s)} \right] + \frac{(1-s)k(s,t)}{t(1-s)+[k(s,t)-s](1-t)}. \label{bound2}
\end{equation}
Below this curve, the decision maker does not accept any object. It follows that in this region the value function at a given point is equal to the time when 
the limiting process intersects the region in which relatively best objects from ordering 2 become acceptable. It follows that $v(s,t) = v\left( c(s,t),t+\frac{[c(s,t)-t][1-t]}{1-s}\right)$, where 
$c(s,t)$ satisfies 
\[
v\left( c(s,t),t+\frac{[c(s,t)-t][1-t]}{1-s}\right) = t+\frac{[c(s,t)-t][1-t]}{1-s}.
\] 

Since the upper boundary of this region is given by an implicit function, numerical estimation of the general form of the 
function is complex. However, it is relatively simple to estimate $v(s,t)$ for specified values of $s$ and $t$. In particular, the optimal expected reward is given by $v(0,0)$. In the limit, the 
trajectory of the process moves along the line $s=t$. Moving along such a trajectory, the DM 
should accept the first candidate to appear when $s\geq 0.5$. It follows that the 
optimal expected reward is $v(0.5,0.5)=0.5$. 

In addition, note that any relatively best object of type 2 should always be accepted when $t\geq t^{\ast}(0)$, even if no 
objects of type 1 have be observed. We can estimate $t^{\ast}(0)$,  using 
$t^{\ast}(0)=v(0,t^{\ast}(0))$.  From Equation (\ref{reg2}), it follows that 
\begin{equation}
1 = \ln \left\{ \frac{t^{\ast}(0)+[1-t^{\ast}(0)] k(0,t^{\ast}(0))} {t^{\ast}(0)} \right\} + \frac{k(0,t^{\ast}(0))}{t^{\ast}(0)+[1-t^{\ast}(0)]k(0,t^{\ast}(0))} .
\label{t0_1}
\end{equation}
Given that no object is chosen, then starting from the state $(0,t^{\ast}(0),0)$ the limiting process will cross into the region where any relatively best object will be accepted in state
$(s^{\ast}(\tilde{t}),\tilde{t},0)$, where  
\begin{equation}
\tilde{t}=t^{\ast}(0)+[1-t^{\ast}(0)]s^{\ast}(\tilde{t})  \label{t0_2}
\end{equation}
and $v(s^{\ast}(\tilde{t}),\tilde{t})=s^{\ast}(\tilde{t})$, which together with Equation (\ref{bound1}) leads to
\begin{equation}
1 = \frac{\tilde{t} (2-s^{\ast}(\tilde{t})-\tilde{t})[\ln (\tilde{t}) - \ln (s^{\ast}\{\tilde{t}\} )]}{\tilde{t}-s^{\ast}(\tilde{t})} . \label{t0_3}
\end{equation}
Equations (\ref{t0_1})-(\ref{t0_3}) were solved numerically using the \textit{dfsane} procedure in the \textit{BB} package written for R (see 
Varadhan and Gilbert \citep{var1}). We obtain $t^{\ast}(0)\approx 0.5082004$ and starting from the state $(0,t^{\ast}(0),0)$, the 
limiting process crosses into the region where  any relatively best object will be accepted at the point $(s^{\ast}(\tilde{t}),\tilde{t})\approx (0.4633848,0.7360929)$.

It should be noted that solving the discrete time problem with $N\geq$ 10 000 gives (for the purposes of visualization) a reasonable approximation to the asymptotic 
solution of the problem. For example, suppose $N=$10 000. If no objects of type 1 have been seen , then a relatively best object of type 2 should be accepted if and only if
at least 5 083 objects of type 2 have appeared. Fig. 1 illustrates the asymptotic form of the solution, estimated on this basis.

\begin{figure}[h]
\centering
\includegraphics[width=.6\textwidth]{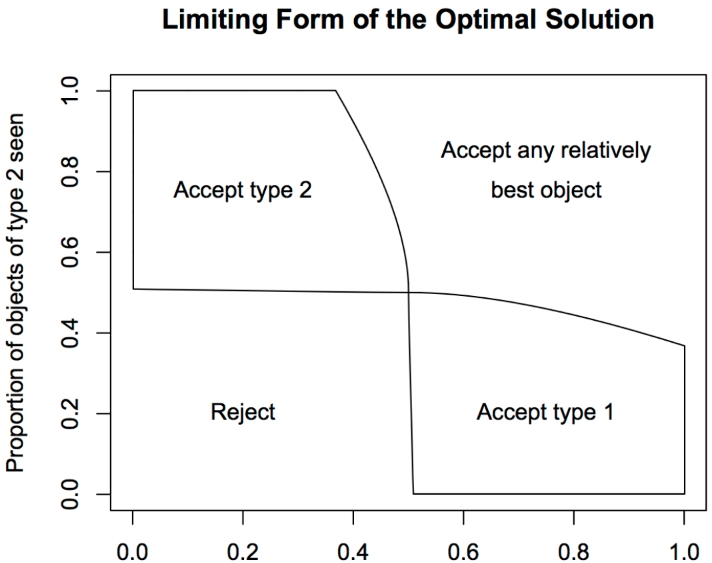}
\label{Fig2}
\caption{Asymptotic Form of the Optimal Strategy}
\end{figure}

\section{Conclusion}

This paper has considered a secretary problem, a stopping problem, in which there are an equal number of 
objects of two different types. Objects of a given type can be fully ranked, but objects of different types cannot be compared. The goal of the DM is choose an object which is the absolute best of its type. Any particular relaization of such a problem can be solved using dynamic programming. Some stopping problems can be solved using the concept of a one-step look-ahead rule (OSLAR). Under such a rule, a DM stops if and only if the expected reward from stopping at the current decision point is at least as great as the expected 
reward from stopping at the next decision point. It was argued that the optimal solution to the problem considered is not an OSLAR. However, given that the process has not yet been stopped, eventually it passes into a region where it is optimal to apply such a rule. 

The form of the optimal strategy was then derived by combining the approaches of dynamic programming and the theory of OSLAR. A particularly strong result on whether or not a candidate should be 
accepted in state $(i,i,1)$ was given. This was followed by a number of results on the asymptotical solution when the number of objects of both types tends to infinity.

There are a number of other directions for future research. One could relax the 
assumption that the number of objects of both types are equal. Suppose there are $N_i$ objects of type $i$, $i=1,2$. One could then 
derive asymptotic results for the case where $N_2 = kN_1$ and $N_1, N_2
\rightarrow \infty$. One could also increase the number of types of object available.
This would increase the complexity of the problem. However, it is likely that the methods applied here could be also used, or adapted, in such problems. One could also
look at similar problems in which the DM has a different goal (e.g. to choose one of the two best objects of a type).





\newpage
\bibliographystyle{plain}
\bibliography{twotwigs_spa}

\newpage
\appendix
\markboth{\shortauthors}{APPENDIX}
\section{Proof of Lemma 5}
 
Let $0\leq a\leq b \leq c \leq d\leq N$, such that $a+d=b+c=n$ and $(a,d)\neq (b,c)$. Then 
$w_{a,d}\leq w_{b,c}.$

\vspace{.1in}

\begin{Proof}
Note that from the definition of $a,b,c,d$ then $a\leq d-2$. 
We shall use recursion with respect to $n$.
 From the symmetry of the problem
$w_{N-1,N}=w_{N,N-1}.$ The largest relevant value of $n$ is $2N-2$.
Assume that $a=N-2,d=N$. In this case we have
\begin{eqnarray*}
w_{N-1,N-1} & = & \frac{1}{N}+\frac{N-1}{N}w_{N-1,N} = \frac{1}{N} +\frac{N-1}{N}\cdot\frac{1}{N} \\
& \geq & \frac{1}{N}+\frac{N-2}{N-1}\cdot\frac{1}{N}=w_{N-2,N}. 
\end{eqnarray*}
Hence our claim is true for $n=2N-2$.

For the rest of the proof, we first show that for $0\leq a\leq d-2 \leq N-2$
\begin{equation}
\frac{N- a - 1}{2N\! - \! a\! - \! d}\cdot\frac{a\! +\! 1}{a\! +\! 2}+\frac{N-d+1}{2N\! -\! a\! -\! d}\cdot\frac{d\! -\! 1}{d}
\geq \frac{N-a}{2N\! -\! a\! - \! d}\cdot\frac{a}{a\! + \! 1}+\frac{N-d}{2N\! -\! a\! -\! d}\cdot\frac{d}{d\! +\! 1}. \label{lm51}
\end{equation}
Let
\[
f(a,d)=\frac{N- a - 1}{2N\! - \! a\! - \! d}\cdot\frac{a\! +\! 1}{a\! +\! 2}+\frac{N-d+1}{2N\! -\! a\! -\! d}\cdot\frac{d\! -\! 1}{d}
- \frac{N-a}{2N\! -\! a\! - \! d}\cdot\frac{a}{a\! + \! 1}-\frac{N-d}{2N\! -\! a\! -\! d}\cdot\frac{d}{d\! +\! 1}.
\]
We have
\begin{eqnarray*}
f(a,\!d) & \! \! = \! \! & \frac{N-a}{2N\! \! -\! a\! -\! d}\! \left[ \! \frac{a\! +\! 1}{a\! +\! 2}-\frac{a}{a\! +\! 1} \right]
\! \! +\! \frac{N-d}{2N\! \! -\! a\! -\! d}\! \left[ \! \frac{d\! -\! 1}{d}\! -\! \frac{d}{d\! +\! 1} \right] \! \! + \! \frac{1}{2N\! \! -\! a\! -\! d}\left[ \frac{d\! -\! 1}{d}-\frac{a\! +\! 1}{a\! +\! 2} \right] \\
& \! \! \geq \! \! & \frac{N-a}{2N-a-d} \cdot \frac{1}{(a+1)(a+2)} - \frac{N-d}{2N-a-d} \cdot 
\frac{1}{d(d+1)} \geq 0.
\end{eqnarray*}
 
Assume now that the recursion hypothesis is true for some $n=a+d+1$, where \newline $1\leq n\leq 2N-2$, and $a\leq d-3$.
Using the recursion hypothesis, from Equation (\ref{recurs2}) we obtain
\begin{eqnarray*}
w_{a+1,d-1} & = & \frac{1}{N}+\frac{N-a-1}{2N-a-d}\cdot\frac{(a+1)w_{a+2,d-1}}{a+2}+\frac{N-d+1}{2N-a-d}\cdot\frac{(d-1)w_{a+1,d} }{d}\\
& \geq & \frac{1}{N}+\frac{N-a-1}{2N-a-d}\cdot\frac{(a+1)w_{a+1,d}}{a+2}+
\frac{N-d+1}{2N-a-d}\cdot\frac{(d-1)w_{a+1,d}}{d} \\
& = & \frac{1}{N}+\left(\frac{N-a-1}{2N-a-d}\cdot\frac{a+1}{a+2}+\frac{N-d+1}{2N-a-d}\cdot\frac{d-1}{d}\right) w_{a+1,d} \\
&  \geq & \frac{1}{N}+\left(\frac{N-a}{2N-a-d}\cdot\frac{a}{a+1}+\frac{N-d}{2N-a-d}\cdot\frac{d}{d+1}\right) w_{a+1,d}\\
& \geq & \frac{1}{N}+\frac{N-a}{2N-a-d}\cdot\frac{aw_{a+1,d}}{a+1}+
\frac{N-d}{2N-a-d}\cdot\frac{dw_{a,d+1}}{d+1}=w_{a,d}. 
\end{eqnarray*}
The first and third inequalities result from the induction hypothesis. The second inequality follows 
from Inequality (\ref{lm51}).

The final case we need to consider is when $d=a+2$ and $a+d\leq N-3$.  
We have 
\begin{eqnarray*}
w_{a,a+2} & = & \frac{1}{N}+\frac{N-a}{2N-2a-2}\cdot\frac{aw_{a+1,a+2}}{a+1}+
\frac{N-a-2}{2N-2a-2}\cdot\frac{(a+2)w_{a,a+3}}{a+3} \\ 
& \leq &
\frac{1}{N}+\frac{N-a}{2N-2a-2}\cdot\frac{aw_{a+1,a+2}}{a+1}+
\frac{N-a-2}{2N-2a-2}\cdot\frac{(a+2)w_{a+1,a+2}}{a+3} \\ 
& = &
\frac{1}{N}+\left[ \frac{N-a}{2N-2a-2}\cdot\frac{a}{a+1}+
\frac{N-a-2}{2N-2a-2}\cdot\frac{(a+2)}{a+3}\right] w_{a+1,a+2} \\ 
& \leq & \frac{1}{N}+\frac{(a+1)w_{a+1,a+2}}{a+2}=w_{a+1,a+1}. 
\end{eqnarray*}
The first equation results from the induction hypothesis. The second inequality results from the fact that when $d=a+2$, Inequality (\ref{lm51}) reduces to 
\[
\frac{a+1}{a+2} \geq \frac{N-a}{2N - 2a - 2}\cdot\frac{a}{a + 1}+\frac{N-a-2}{2N - 2a - 2}\cdot\frac{a+2}{a + 3}.
\]
\end{Proof}

\section{Proof of Lemma 6}

Suppose that $j\leq i<N$. It follows that
\[
\frac{iw_{i+1,j}}{i+1} \geq \frac{jw_{i,j+1}}{j+1}.
\]

\begin{Proof} First note that by symmetry $w_{i+1,i}=w_{i,i+1}$ and thus when $i=j$ the inequality in the lemma is 
satisfied with equality. Otherwise,
the proof is by recursion according to the total number of objects that have appeared, $n$, for $n\leq 2N-3$. 

When $n=2N-3$, the only case that we need to consider is $i=N-1, j=N-2$. It is necessary to show that
\begin{equation}
\frac{(N-1)w_{N,N-2}}{N} \geq \frac{(N-2)w_{N-1,N-1}}{N-1}. \label{ineq}
\end{equation}
 Using Equation (\ref{recurs2}), we obtain
\begin{eqnarray*}
w_{N,N-2} & = & \frac{1}{N} + \frac{(N-2)w_{N,N-1}}{N-1} = \frac{2N-3}{N(N-1)} \\
w_{N-1,N-1} & = & \frac{1}{N} + \frac{(N-1)w_{N,N-1}}{N} = \frac{2N-1}{N^2}. 
\end{eqnarray*}
Hence,
\begin{eqnarray*}
\frac{(N-2)w_{N-1,N-1}}{N-1} & = & \frac{(N-2)(2N-1)}{N^2 (N-1)} = \frac{2N^2-5N+2}{N^2(N-1)} \\
\frac{(N-1)w_{N,N-2}}{N} & = & \frac{(2N-3)(N-1)}{N^2 (N-1)} = \frac{2N^2-5N+3}{N^2(N-1)}.
\end{eqnarray*}
It is thus clear that Inequality (\ref{ineq}) holds for any $N\geq 2$ (which is necessary for 
there to be three objects yet to appear).

Now suppose that the lemma holds for $n\geq 2N-k+1$. We now show that the lemma holds 
for $n=2N-k$. The lemma then follows by recursion.

Suppose that $n=i+j=2N-k$ and $i>j$. Define $\delta = \frac{1}{n+1}$. From Equation (\ref{recurs2}) 
\begin{eqnarray*}
w_{i+1,j} & = & \frac{1}{N} + \frac{(i+1)\alpha_{i+1,j} w_{i+2,j}}{i+2} + \frac{j (1-\alpha_{i+1,j} )w_{i+1,j+1}}{j+1} \\
& \geq & \frac{1}{N} + \frac{jw_{i+1,j+1}}{j+1},
\end{eqnarray*}
where $0\leq \alpha_{i+1,j}=\frac{N-i-1}{2N-i-j-1}<0.5$ and the inequality results from the induction assumption. Similarly,
\begin{eqnarray*}
w_{i,j+1} & = & \frac{1}{N} + \frac{i(\alpha_{i,j} +\delta)w_{i+1,j+1}}{i+1} + \frac{(j+1)(1-\alpha_{i,j} - \delta )w_{i,j+2}}{j+2} \\
& \leq & \frac{1}{N} + \frac{iw_{i+1,j+1}}{i+1}.
\end{eqnarray*}
Multiplying the above inequalities by $\frac{i}{i+1}$ and $\frac{j}{j+1}$, respectively, it follows that
\begin{eqnarray*} 
\frac{iw_{i+1,j}}{i+1} - \frac{jw_{i,j+1}}{j+1} \geq \frac{1}{N}\left[ \frac{i}{i+1}-\frac{j}{j+1}\right] \geq 0.
\end{eqnarray*}
The final inequality results from the fact that $i>j$.
\end{Proof}

\section{Proof of Lemma \ref{1N}}

For $1\leq i,j \leq N-1$, $|v_{i+1,j}-v_{i,j+1}| \leq \frac{1}{N}$.

\begin{Proof} The proof is by recursion on the number of offers that have appeared, $n=i+j$. 

For $n=2N-1$, $v_{N, N-1}=v_{N-1,N}=\frac{1}{N}$. Hence, the hypothesis is true for $n=2N-1$. 

Assume that this lemma is true for $n = 2N-s$, where $s\geq 1$. 
Given that $i$ objects of type 1 and $j$ objects of type 2 have appeared, let 
$v_{i,j|k}$ be the optimal reward given that the next object to appear is of 
type $k$. 

\begin{eqnarray} 
v_{i-1,j} & = & \alpha_{i-1,j} v_{i-1,j|1} + (1-\alpha_{i-1,j})v_{i-1,j|2} \label{eq4f1} \\
v_{i,j-1} \! \! & = & \! \! \left[ \alpha_{i-1,j} -\frac{1}{N-n+1}\right] \! v_{i,j-1|1} \! +\left[1- \alpha_{i-1,j} \! +\! \frac{1}{N-n+1}\right] \! v_{i,j-1|2}.
\label{eq4f2}
\end{eqnarray}

Note that
\begin{eqnarray*}
v_{i-1,j|1} & = & \max\left\{ v_{i,j}, \frac{1}{N}+\frac{(i-1)v_{i,j}}{i}\right\} \\
v_{i-1,j|2} & = & \max\left\{ v_{i-1,j+1}, \frac{1}{N}+\frac{jv_{i-1,j+1}}{j+1} \right\} .
\end{eqnarray*}

Firstly, suppose that $v_{i,j} \geq \frac{1}{N}+\frac{(i-1)v_{i,j}}{i}$ and  $v_{i+1,j-1} \geq \frac{1}{N}+\frac{iv_{i+1,j-1}}{i+1}$. It follows that 
\[
|v_{i-1,j|1}-v_{i,j-1|1}| = |v_{i,j}-v_{i+1,j-1}| \leq \frac{1}{N},
\]
where the inequality results from the recursion hypothesis.

Secondly, suppose that  $v_{i,j} < \frac{1}{N}+\frac{(i-1)v_{i,j}}{i}$ and  $v_{i+1,j-1} < \frac{1}{N}+\frac{iv_{i+1,j-1}}{i+1}$. It follows that $v_{i+1,j-1} \leq \frac{i+1}{N}$. We have
\begin{eqnarray*}
|v_{i-1,j|1}-v_{i,j-1|1}| & = & \left| \frac{(i-1)v_{i,j}}{i}-\frac{iv_{i+1,j-1}}{i+1}\right|  \\
& \leq & \frac{i-1}{i} \left| v_{i,j}-v_{i+1,j-1}\right| + \frac{v_{i+1,j-1}}{i(i+1)} \\
& \leq & \frac{i-1}{iN} + \frac{1}{i(i+1)} \times \frac{i+1}{N} \leq \frac{1}{N}.
\end{eqnarray*}
Thirdly, suppose that $v_{i,j} \geq \frac{1}{N}+\frac{(i-1)v_{i,j}}{i}$ and  $v_{i+1,j-1} < \frac{1}{N}+\frac{iv_{i+1,j-1}}{i+1}$. In this case, since $v_{i+1,j-1}< \frac{i+1}{N}$ and $v_{i,j}\geq \frac{i}{N}$, it follows from the recursion hypothesis that $\frac{i}{N}\leq v_{i,j}\leq \frac{i+2}{N}$. We have 
\begin{eqnarray*}
|v_{i-1,j|1}-v_{i,j-1|1}| & = & \left|v_{i,j}-\frac{iv_{i+1,j-1}}{i+1}-\frac{1}{N}\right|  \\
& \leq & \frac{i}{i+1} \left| v_{i,j}-v_{i+1,j-1}\right| + \left|\frac{v_{i,j}}{i+1} - \frac{1}{N} \right| .\\
\end{eqnarray*}
Note that 
\[
 \frac{-1}{(i+1)N} = \frac{i}{(i+1)N} - \frac{1}{N} \leq \frac{v_{i,j}}{i+1} - \frac{1}{N} \leq \frac{i+2}{(i+1)N} - \frac{1}{N} = \frac{1}{(i+1)N}.
\]
Hence, 
\[
|v_{i-1,j|1}-v_{i,j-1|1}|  \leq  \frac{i}{(i+1)N} +   \frac{1}{(i+1)N} = \frac{1}{N}.
\]
Finally, note that the combined conditions $v_{i,j} < \frac{1}{N} + \frac{(i-1)v_{i,j}}{i}$ and  $v_{i + 1,j - 1} \geq \frac{1}{N} + \frac{iv_{i+1,j-1}}{i+1}$ contradict the induction hypothesis, since it follows that 
 $v_{i + 1,j - 1} \geq \frac{i+1}{N}$ and $v_{i,j} < \frac{i}{N}$. 
Hence, $|v_{i-1,j|1}-v_{i,j-1|1}|\leq \frac{1}{N}$.

Note that by symmetry $v_{i,j|2}=v_{j,i|1}$. Hence, it follows from the above arguments that 
$|v_{i-1,j|2}-v_{i,j-1|2}| \leq \frac{1}{N}$.

Let $\delta = \alpha_{i-1,j}-\alpha_{i,j-1}=\frac{1}{N-n+1}$.
From Equations (\ref{eq4f1}) and (\ref{eq4f2}), it follows that
\begin{eqnarray}
|v_{i - 1,j} -  v_{i,j - 1}|  & = & \left| \alpha_{i-1,j} v_{i-1,j|1} + (1 - \alpha_{i-1,j})v_{i-1,j|2}+ \right. \nonumber \\
& & - \left. ( \alpha_{i-1,j} - \delta ) v_{i,j-1|1} - (1 - \alpha_{i-1,j} + \delta ) v_{i,j-1|2} \right| \nonumber \\ 
 & \leq & ( \alpha_{i-1,j} - \delta ) \left| v_{i-1,j|1} -  v_{i,j-1|1}\right| + (1 - \alpha_{i-1,j} ) \left| v_{i - 1,j|2} - v_{i,j - 1|2}\right| + \nonumber \\
& & +
\delta \left| v_{i-1,j|1} - v_{i,j - 1|2}\right| \nonumber \\
 & \leq & \frac{1}{N} \left[ 1 -\delta \right] + \delta \left| v_{i-1,j|1}-v_{i,j-1|2}\right| .\label{eq4f3}
\end{eqnarray}
Suppose that in both of the states $(i,j,1)$ and $(i,j,2)$, the DM should reject a candidate. It follows that 
$v_{i-1,j|1}=v_{i,j-1|2}=v_{i,j}$. 

Now suppose that in both of the states $(i,j,1)$ and $(i,j,2)$, the DM should accept a candidate. It follows that
\begin{eqnarray*}
|v_{i-1,j|1}-v_{i,j-1|2}| & = & \left| \left[ \frac{1}{N}+\frac{(i-1)v_{i,j}}{i}\right] - \left[ \frac{1}{N}+\frac{(j-1)v_{i,j}}{j}\right] \right| \\
& & =  v_{i,j} \left| \left( \frac{i-1}{i}-\frac{j-1}{j} \right) \right| 
 = v_{i,j} \left| \frac{i-j}{ij}\right| .
\end{eqnarray*}
From the symmetry between $i$ and $j$, we may assume without loss of generality that $i\geq j$. Since a 
candidate should be accepted in state $(i,j,2)$, 
$v_{i,j} \leq \frac{j}{N}$. Thus
\[
|v_{i-1,j|1}-v_{i,j-1|2}| \leq \frac{(i-j)j}{ijN} = \frac{i-j}{iN} \leq \frac{1}{N}.
\]
Assume that a DM should accept a candidate in state $(i,j,1)$, but not in state $(i,j,2)$. Hence, $i>j$ and 
$\frac{j}{N} < v_{i,j} \leq \frac{i}{N}$.
\[
|v_{i-1,j|1}-v_{i,j-1|2}| = \left| \frac{1}{N}+\frac{(i-1)v_{i,j}}{i}- v_{i,j}\right| = \left| \frac{1}{N} - \frac{v_{i,j}}{i}\right| \leq \frac{1}{N}.
\]
It follows from an analogous argument that when a DM should accept a candidate in state $(i,j,2)$, but not in state $(i,j,1)$, then 
$|v_{i-1,j|1}-v_{i,j-1|2}| \leq \frac{1}{N}$. Hence, for all $1\leq i,j \leq N$, 
$|v_{i-1,j|1}-v_{i,j-1|2}| \leq \frac{1}{N}$. 

It thus follows from Equation (\ref{eq4f3}) that $|v_{i - 1,j} -  v_{i,j - 1}|  \leq \frac{1}{N}$ and Lemma \ref{1N} follows by recursion.
\end{Proof}

\section{Proof of Lemma \ref{1n}} 

Assume that it is optimal to accept an offer in state 
$(i+1,j,1)$, where $j\leq i<N$. It follows that 
\begin{equation}
\frac{1}{N} + \frac{i v_{i+1,j}}{i+1} \geq v_{i,j+1} \label{L1n}.
\end{equation}

\begin{Proof} 
 First, assume that $(i,j)\in B_3$, i.e. starting from the state $(i,j,0)$ it is always optimal to accept the next candidate to appear. From Lemma \ref{w_dec} it follows that
\[
\frac{1}{N} + \frac{i v_{i+1,j}}{i+1} \geq \frac{1}{N} + \frac{j v_{i,j+1}}{j+1}.
\]
In this case, the hypothesis follows from the optimality criterion
\[
v_{i,j+1} \leq \frac{j+1}{N} \Rightarrow \frac{1}{N} + \frac{j v_{i,j+1}}{j+1} \geq v_{i,j+1}.
\]
Note that from Theorem \ref{acc_next}, when $(i,j)\in B_3$ it follows that $v_{i,j}\geq v_{i,j+1}$.

When $i\geq i^{\dagger}_N$ a candidate should not be accepted in state 
$(i,j^{\circ} (i),2)$, but always be accepted in state 
$(i, j^{\circ} (i)+1,k)$, where $k\in \{ 1,2\}$. 
Now we consider starting from a state $(i,j,0)$, $j<i^{\dagger}_N \leq i<N$, that borders the set $B_3$, i.e. $(i+1,j)$ and/or $(i,j+1)\in B_3$, but $(i,j)\notin B_3$.
Consider a state of the form $(i, j,0)$ such that both $(i+1,j)$ and $(i,j+1)$ are in $B_3$.  It follows that a candidate should be accepted in both of the 
states $(i+1,j+1,1)$ and $(i,j+2,2)$. Note that by assumption a candidate should be accepted in state $(i+1,j,1)$. We have
\begin{eqnarray}
\frac{1}{N} \! + \! \frac{iv_{i+1,j}}{i+1}\! -\! v_{i,j+1}
\! & \geq & \! \frac{1}{N} \! + \! \frac{i}{i\! +\! 1} \left[ \alpha_{i+1,j} \! \left( \frac{1}{N} \! + \!
\frac{(i\! +\! 1)v_{i+2,j}}{i+2}\right) \! + \! (1\! -\! \alpha_{i+1,j})v_{i+1,j+1}\right] \! + \! \nonumber \\
& & -(\alpha_{i+1,j} \! +\! \epsilon ) \! \left[ \! \frac{1}{N} \! +\! \frac{iv_{i+1,j+1}}{i+1} \right] 
\! - \! (1\! -\! \alpha_{i+1,j}\! -\! \epsilon ) \! \left[ \! \frac{1}{N} \! + \! \frac{(j\! +\! 1)v_{i,j+2}}{j+2} 
\right] \nonumber \\
& = & \frac{i\alpha_{i+1,j}}{i+1}\left[ \frac{1}{N} +\frac{(i+1)v_{i+2,j}}{i+2}
-v_{i+1,j+1}\right] + \nonumber \\ 
& & + (1-\alpha_{i+1,j}-\epsilon )\left[ \frac{iv_{i+1,j+1}}{i+1}
-\frac{(j+1)v_{i,j+2}}{j+2}\right] \geq 0, \label{border1}
\end{eqnarray}
where $\epsilon = \frac{1}{2N-i-j-1}$. The first term in square brackets in the final inequality 
is non-negative because $(i+1,j)\in B_3$. 
The fact that the second term in square brackets is non-negative results from Lemma 6 and the facts that $i>j$ and $(i,j+1)\in B_3$. 
Note that for such states
\[
v_{i,j}= \alpha_{i,j} \left[ \frac{1}{N}+\frac{iv_{i+1,j}}{i+1}\right] 
+(1-\alpha_{i,j})v_{i,j+1}.
\]
Hence, it follows that $v_{i,j}\geq v_{i,j+1}$, since $v_{i,j}$ is the weighted average of 
$v_{i,j+1}$ and another term which is at least as large as $v_{i,j+1}$. 

Let $t(j)$ be the largest
integer $i$ for a given $j$, where $j^{\circ}(N)\leq j\leq j^{\circ}(i^{\bullet}_N)$, such that $(i,j)\notin B_3$. From the argument above, the lemma holds for $(i,j)=(t(j),j)$.  
Now we consider starting in a state of the form $(i, j,0)$ such that $j\geq j^{\circ}(N)$, a candidate should be accepted in state $(i+1,j,1)$, $(i+1,j)\notin B_3$, but $(i,j+1)\in B_3$. It follows that a candidate should be 
accepted in both of the states $(i+1,j+1,1)$ and $(i,j+2,2)$. Hence, for a given $j$, the set of $(i,j)$ satisfying this condition is of the form  $\{ (s(j),j), (s(j)+1,j),\ldots , (t(j)-1,j)\}$.  Adapting the derivation 
of Inequality (\ref{border1}) for $(i,j)=(t(j)-1,j)$, it follows that the lemma holds for
state $(t(j)-1,j)$. The only change in the argument is that the first term in square brackets in Inequality (\ref{border1}) is non-negative due to the lemma holding 
for $(i,j)=(t(j),j)$.
Arguing analogously, it follows by recursion that the lemma holds 
for each state in the set $(s(j),j), (s(j)+1,j),\ldots , (t(j)-1,j)$. 
Hence, the lemma holds for any
$(i,j)\notin B_3$ such that $(i,j+1)\in B_3$.

Assume that $i> j$, it is impossible for $(i+1,j)\in B_3$ and $(i,j+1)\notin B_3$. The first condition 
requires that $v_{i+1,j+1}\leq \frac{j+1}{N}$. From Lemma \ref{1N} it follows
that $v_{i,j+2}\leq \frac{j+2}{N}$. Hence, starting from state $(i,j+1,0)$, any  candidate which appears at the next moment, i.e. in state $(i+1,j+1,1)$ or 
$(i,j+2,2)$ should be accepted.

From the above arguments, it follows that 
Inequality (\ref{L1n}) holds whenever $(i+1,j)\in B_3$ and/or $(i,j+1)\in B_3$.
Let $U_0$ be the set of $(i,j)$, such that $i> j$, a candidate should be accepted in state $(i,j,1)$, $(i,j)\notin B_3$ but 
$(i+1,j)\in B_3$ and/or $(i,j+1)\in B_3$. 

Note that $v_{i,j} \geq v_{i,j+1}$ for $(i,j)\in U_0$. Hence, if a candidate 
should be accepted in state $(i,j,1)$, i.e. $v_{i,j}\leq \frac{i}{N}$, it follows that
a candidate should be accepted in state $(i,j+1,1)$, since 
$v_{i,j+1}\leq \frac{i}{N}$. From Lemma \ref{1N}, $v_{i+1,j}\leq \frac{i+1}{N}$ and thus a candidate should be accepted in state $(i+1,j,1)$.

Finally, we consider starting from a state $(i,j,0)$ such that $i>j$, a candidate should be 
accepted in state $(i+1,j,1)$, but neither $(i+1,j)$ nor $(i,j+1)$ belong to $B_3$.
Since a candidate should be accepted in state $(i+1,j+1,1)$
\begin{eqnarray}
\frac{1}{N} \! + \! \frac{iv_{i+1,j}}{i+1}\! -\! v_{i,j+1}
\! & \geq & \! \frac{1}{N} \! + \! \frac{i}{i\! +\! 1} \left[ \alpha_{i+1,j} \! \left( \frac{1}{N} \! + \!
\frac{(i\! +\! 1)v_{i+2,j}}{i+2}\right) \! + \! (1\! -\! \alpha_{i+1,j})v_{i+1,j+1}\right] \! + \!  \nonumber \\
& & -(\alpha_{i+1,j} \! +\! \epsilon ) \! \left[ \! \frac{1}{N} \! +\! \frac{iv_{i+1,j+1}}{i+1} \right] 
\! - \! (1\! -\! \alpha_{i+1,j}\! -\! \epsilon ) v_{i,j+2}
 \nonumber \\
& = & \frac{i\alpha_{i+1,j}}{i+1}\left[ \frac{1}{N} +\frac{(i+1)v_{i+2,j}}{i+2}
-v_{i+1,j+1}\right] + \nonumber \\ 
& & + (1-\alpha_{i+1,j}-\epsilon )\left[ \frac{1}{N}+ \frac{iv_{i+1,j+1}}{i+1}
-v_{i,j+2} \right]. \label{border2}
\end{eqnarray}
Define the set $U_m$ iteratively for $m\geq 1$ as
\[
U_m = U_{m-1} \cup \{ (i,j): (i+1,j)\in U_{m-1}, (i,j+1)\in U_{m-1}, 
v_{i+1,j}\leq \frac{i+1}{N} \} . 
\] 
It follows that $U_m$ denotes the set of $(i,j)$ such that in $m$ steps of the search 
process the state is guaranteed to have reached (or passed through) the set $U_0$ 
and it is optimal to accept a candidate in state $(i+1,j,1)$.
From Inequality
(\ref{L1n}), for $(i,j)\in U_0$, we have $v_{i+1,j}\geq v_{i+1,j+1}$. Hence, given that a candidate should be accepted in state $(i+1,j,1)$, it is also optimal to accept a candidate in state $(i+1,j+1,1)$. From Inequality (\ref{border2}) it then follows that Inequality
(\ref{L1n}) holds for $(i,j)\in U_1$. Similarly, if this inequality holds for $(i,j)\in U_m$, then 
it holds for $(i,j)\in U_{m+1}$. It then follows by induction that 
this inequality holds whenever it is optimal to accept a candidate in state $(i+1,j,1)$.
In addition,  whenever $v_{i,j} \leq \frac{i}{N}$, then $v_{i+1,j} \leq \frac{i+1}{N}$ and from the dynamic programming equation it follows that  $v_{i,j}\geq v_{i,j+1}$ for such $(i,j)$. Hence, if it is optimal to accept a candidate in state 
$(i,j,1)$, then it is optimal to accept a candidate in states $(i+1,j,1)$ and 
$(i,j+1,1)$. Thus it follows by induction that if the search process enters a region in 
which it is optimal to accept a candidate of type 1, it will be optimal to accept a 
candidate of type 1 at any later moment. By symmetry, an analogous result holds for 
candidates of type 2.
\end{Proof}

\section{Proof of Lemma \ref{npol1}}
For all $N$, $i^{\bullet}_N \geq \lceil \frac{N+1}{2}\rceil$.

For $N=2m$ and $N=2m+1$ we shall prove that it is not optimal to stop before time $\left\lceil \frac{N+1}{2}\right\rceil$. 

First we assume that $N=2m$, i.e. $N$ is even.
We shall prove that in the state $(m,m,1)$ it is not optimal to stop. Namely we are going to show that
\begin{equation}\label{before1}
w_{m,m}>\frac{m}{2m}=\frac{1}{2}.
\end{equation}
We have
$$w_{m,m}=2\frac{m^2}{4m^2}\sum_{k=m}^{2m-1}\frac{1}{k}+\frac{m^2}{4m^2}\sum_{u=1}^{m}\sum_{v=1}^{m}\frac{{2m-u-1\choose m-1}{2m-v-1\choose m-1}}{{2m-1\choose m-1}{2m-1\choose m-1}}\frac{2}{u+v}.$$
Thus inequality (\ref{before1}) takes the form
\begin{equation}\label{before2}
\sum_{k=m}^{2m-1}\frac{1}{k}+\sum_{u=1}^{m}\sum_{v=1}^{m}\frac{{2m-u-1\choose m-1}{2m-v-1\choose m-1}}{{2m-1\choose m-1}^2}\frac{1}{u+v}>1.
\end{equation}
 Let us note that 
$$\sum_{u=1}^{m}{2m-u-1\choose m-1}={2m-1\choose m-1}.$$

Thus if we set
$$\alpha(m,u):=\frac{{2m-u-1\choose m-1}}{{2m-1\choose m-1}},$$
and
$$\alpha(m,v):=\frac{{2m-v-1\choose m-1}}{{2m-1\choose m-1}},$$
then
$$\sum_{u=1}^{m}\alpha(m,u)=1\;\;\textrm{and}\;\;\sum_{v=1}^{m}\alpha(m,v)=1$$
and
$$\sum_{u,v=1}^{m}\alpha(m,u)\alpha(m,v)=1.$$
Thus inequality (\ref{before2}) can be written as
\begin{equation}\label{before3}
\sum_{k=m}^{2m-1}\frac{1}{k}+\sum_{u,v=1}^{m}\alpha(m,u)\alpha(m,v)\frac{1}{u+v}>\sum_{u,v=1}^{m}\alpha(m,u)\alpha(m,v).
\end{equation}
We can still transform (\ref{before3}) to the form
\begin{equation}\label{before4}
\sum_{k=m}^{2m-1}\frac{1}{k}>\sum_{u,v=1}^{m}\alpha(m,u)\alpha(m,v)\left(1-\frac{1}{u+v}\right).
\end{equation}
The right-hand side of (\ref{before4}) can be expressed as
$$\sum_{k=m}^{2m-1}\frac{1}{k}>\sum_{u,v=1}^{m}\alpha(m,u)\alpha(m,v)\left(1-\frac{1}{u+v}\right)=
\sum_{k=2}^{2m}\;\sum_{1\leq u,v\leq m, u+v=k} \! \! \! \! \! \! \! \! \! \! \! \! \! \! \alpha(m,u)\alpha(m,v)\left(1-\frac{1}{k}\right).$$
Because there is no more than $k-1$ pairs $(u,v)$ such that $u+v=k$ we obtain for the right-hand side of  (\ref{before4}) we have the following estimation
$$\sum_{k=2}^{2m}\;\sum_{1\leq u,v\leq m, u+v=k}\alpha(m,u)\alpha(m,v)\left(1-\frac{1}{k}\right)<$$
$$\left(\frac{m}{2m-1}\right)^2\left(1-\frac{1}{2}\right)+2\left(\frac{m}{2m-1}\right)^2\cdot\frac{1}{2}\left(1-\frac{1}{3}\right)+
\sum_{k=4}^\infty\left(\frac{1}{2}\right)^k\left(\frac{k-1}{k}\right)(k-1)=$$
$$\sum_{k=2}^\infty\left(\frac{1}{2}\right)^k\left(\frac{k-1}{k}\right)(k-1)+\left[\left(\frac{m}{2m-1}\right)^2\left(1-\frac{1}{2}\right)-\left(\frac{1}{2}\right)^2\left(1-\frac{1}{2}\right)\right]+$$
$$\left[2\left(\frac{m}{2m-1}\right)^2\cdot\frac{1}{2}\left(1-\frac{1}{3}\right)-2\left(\frac{1}{2}\right)^3\left(1-\frac{1}{3}\right)\right]<$$
$$\sum_{k=2}^\infty\left(\frac{1}{2}\right)^k\left(\frac{k-1}{k}\right)(k-1)+\frac{1}{2(2m-1)^2}+\frac{1}{3(2m-1)^2}=$$
$$\sum_{k=2}^\infty\left(\frac{1}{2}\right)^k\left(\frac{k-1}{k}\right)(k-1)+\frac{5}{6(2m-1)^2}.$$
As 
$$\sum_{k=2}^\infty\left(\frac{1}{2}\right)^k\left(\frac{k-1}{k}\right)(k-1)=1-\ln2,$$
it is enough to notice now that
$$\sum_{k=1}^{2m-2}\frac{1}{k}+\left(\frac{1}{2m-1}-\frac{5}{6(2m-1)^2}\right)>\ln 2,$$
and this follows by the well known estimation
$$\sum_{k=1}^{2m-1}\frac{1}{k}-\ln 2>\frac{1}{4m}.$$

\vspace{0.2 cm}

For $N=2m+1$ the proof that $w_{m,m}>\frac{m}{2m+1}$ is analogous.

\section{Proof of Lemma \ref{npol2}}

We first consider that case that $N$ is an even integer such that $N=2M\geq 4$. Let $i,j\geq 1$. Denote
$$Z_1 (i,j)=\frac{i}{N^2}\sum_{r=0}^{N-j-1}\frac{j}{j+r},$$
$$Z_2 (i,j)=\frac{j}{N^2}\sum_{r=0}^{N-i-1}\frac{i}{i+r},$$
$$Z_3 (i,j)=\frac{ij}{N^2}\sum_{u=2}^{N-i+1}\sum_{u=2}^{N-j+1}\frac{{N-u\choose i-1}}{{N-1\choose i-1}}\frac{{N-v\choose j-1}}{{N-1\choose j-1}}\frac{2}{u+v-2}.$$

It follows from Equation \ref{OSLAR} that $w_{i,j}=Z_1 (i,j)+Z_2 (i,j)+Z_3 (i,j)$.

We have
$$Z_1 \left( \frac{N}{2}+1,\frac{N}{2}+1 \right) =Z_2 \left( \frac{N}{2}+1,\frac{N}{2}+1 \right)=\left(\frac{1+\frac{N}{2}}{N}\right)^2\sum_{s=\frac{N}{2}+1
}^{N-1}\frac{1}{s}$$
$$\leq \left(\frac{1+\frac{N}{2}}{N}\right)^2\left(\ln(N-1)-\ln\left(\frac{N}{2}\right)\right).$$

We have
$$\frac{{N-u\choose \frac{N}{2}}}{{N-1\choose \frac{N}{2}}}=\frac{\left(\frac{N}{2}-1\right)\left(\frac{N}{2}-2)\right)\ldots\left(\frac{N}{2}-u+1\right)}{(N-1)(N-2)\ldots(N-u+1)}<\left(\frac{\frac{N}{2}-1}{N-1}\right)^{u-1},$$
and
$$\frac{{N-v\choose i-1}}{{N-1\choose i-1}}=\frac{\left(\frac{N}{2}-1\right)\left(\frac{N}{2}-2)\right)\ldots\left(\frac{N}{2}-v+1\right)}{(N-1)(N-2)\ldots(N-u+1)}<\left(\frac{\frac{N}{2}-1}{N-1}\right)^{v-1}.$$
Hence
$$Z_3\left(\frac{N}{2}+1,\frac{N}{2}+1\right)<\left(\frac{1+\frac{N}{2}}{N}\right)^2\sum_{u=2}^{\infty}\sum_{v=2}^{\infty}\left(\frac{\frac{N}{2}-1}{N-1}\right)^{u+v-2}\frac{2}{u+v-2}$$
$$=2\left(\frac{1+\frac{N}{2}}{N}\right)^2\sum_{k=1}^{\infty}\sum_{l=1}^{\infty}\frac{1}{k+l}\left(\frac{\frac{N}{2}-1}{N-1}\right)^{k+l} \! \! \! =
2\left(\frac{1+\frac{N}{2}}{N}\right)^2\sum_{m=2}^{\infty}\frac{m-1}{m}\left(\frac{\frac{N}{2}-1}{N-1}\right)^{m}$$
$$=2\left(\frac{1+\frac{N}{2}}{N}\right)^2\left(\frac{\left(\frac{\frac{N}{2}-1}{N-1}\right)^2}{1-\left(\frac{\frac{N}{2}-1}{N-1}\right)}-\sum_{m=2}^{\infty}\frac{1}{m}\left(\frac{\frac{N}{2}-1}{N-1}\right)^{m}\right)$$
$$=2\left(\frac{1+\frac{N}{2}}{N}\right)^2\frac{\left(\frac{\frac{N}{2}-1}{N-1}\right)^2}{1-\left(\frac{\frac{N}{2}-1}{N-1}\right)}-2\left(\frac{1+\frac{N}{2}}{N}\right)^2\sum_{m=2}^{\infty}\frac{1}{m}\left(\frac{\frac{N}{2}-1}{N-1}\right)^{m}$$
$$=2\left(\frac{1+\frac{N}{2}}{N}\right)^2\frac{\left(\frac{\frac{N}{2}-1}{N-1}\right)^2}{1-\left(\frac{\frac{N}{2}-1}{N-1}\right)}-2\left(\frac{1+\frac{N}{2}}{N}\right)^2\sum_{m=2}^{\infty}\frac{1}{m}x^{m}|_{x=\frac{\frac{N}{2}-1}{N-1}},$$
further we set $a=\frac{\frac{N}{2}-1}{N-1}$,
$$=2\left(\frac{1+\frac{N}{2}}{N}\right)^2\frac{\left(\frac{\frac{N}{2}-1}{N-1}\right)^2}{1-\left(\frac{\frac{N}{2}-1}{N-1}\right)}-2\left(\frac{1+\frac{N}{2}}{N}\right)^2\sum_{m=2}^{\infty}\int_0^a x^{m-1}\,dx$$
$$=2\left(\frac{1+\frac{N}{2}}{N}\right)^2\frac{\left(\frac{\frac{N}{2}-1}{N-1}\right)^2}{1-\left(\frac{\frac{N}{2}-1}{N-1}\right)}-2\left(\frac{1+\frac{N}{2}}{N}\right)^2\int_0^a\sum_{m=1}^{\infty} x^{m}\,dx$$
$$=2\left(\frac{1+\frac{N}{2}}{N}\right)^2\frac{\left(\frac{\frac{N}{2}-1}{N-1}\right)^2}{1-\left(\frac{\frac{N}{2}-1}{N-1}\right)}-2\left(\frac{1+\frac{N}{2}}{N}\right)^2\int_0^a \frac{x}{1-x}\,dx$$
$$=2\left(\frac{1+\frac{N}{2}}{N}\right)^2\frac{\left(\frac{\frac{N}{2}-1}{N-1}\right)^2}{1-\left(\frac{\frac{N}{2}-1}{N-1}\right)}-2\left(\frac{1+\frac{N}{2}}{N}\right)^2\int_0^a (-1)+\frac{1}{1-x}\,dx$$
$$=2\left(\frac{1+\frac{N}{2}}{N}\right)^2\frac{\left(\frac{\frac{N}{2}-1}{N-1}\right)^2}{1-\left(\frac{\frac{N}{2}-1}{N-1}\right)}-2\left(\frac{1+\frac{N}{2}}{N}\right)^2\left(-a-\ln(1-a)\right)$$
$$=\frac{1}{2}\cdot\frac{1\! +\! \frac{N}{2}}{N}\cdot\frac{N^3\! -\! 4N\! -\! 2N^2\! +\! 8}{N^3-N^2}+\frac{1}{2}\cdot\frac{1\! +\! \frac{N}{2}}{N}\cdot\frac{N^2\! -\! 4}{N^2\! -\! N}-\! 2\! \left(\frac{1\! +\! \frac{N}{2}}{N}\right)^2 \! \left(\! \ln(N\! -\! 1)\! -\! \ln \! \left(\frac{N}{2}\right)\! \right)$$
$$=\frac{1+\frac{N}{2}}{N}\cdot\frac{N^3-N^2-4N+4}{N^3-N^2}-2\left(\frac{1+\frac{N}{2}}{N}\right)^2\left(\ln(N-1)-\ln\left(\frac{N}{2}\right)\right)\leq$$
$$\frac{1+\frac{N}{2}}{N}-2\left(\frac{1+\frac{N}{2}}{N}\right)^2\left(\ln(N-1)-\ln\left(\frac{N}{2}\right)\right)$$
Hence, for $M=\frac{N}{2}$
$$w_{M+1,M+1}=2Z_1\left(1+\frac{N}{2},1+\frac{N}{2}\right)+Z_3\left(1+\frac{N}{2},1+\frac{N}{2}\right)\leq\frac{1+M}{N}.$$

Now let us assume that $N$ is odd, $N=2M+1\geq 3$ and $i,j\geq 1$. As in the case where $N$ is even
$$w_{M+1,M+1}=2Z_1(M+1,M+1)+Z_3(M+1,M+1).$$

We have
$$Z_1\left(M+1,M+1\right)=\left(\frac{M+1}{2M+1}\right)^2\sum_{s=M+1}^{2M}\frac{1}{s}.$$
Also,
\[Z_3(M+1,M+1) = 2\cdot\left(\frac{M+1}{2M+1}\right)^{2}\sum_{u=2,v=2}^{M+1}\frac{{2M+1-u\choose M}}{{2M\choose M}}\frac{{2M+1-v\choose M}}{{2M\choose M}}\frac{1}{u+v-2}= \]
\[2\cdot\left(\frac{M+1}{2M+1}\right)^{2}\sum_{u=1,v=1}^{M}\frac{{2M-u\choose M}}{{2M\choose M}}\frac{{2M-v\choose M}}{{2M\choose M}}\frac{1}{u+v}. \]

In addition,
$$\frac{{2M-u\choose M}}{{2M\choose M}}=\frac{M}{2M}\frac{M-1}{2M-1}\ldots\frac{M-u+1}{2M-u+1}= \prod_{i=0}^{u-1}\frac{M-i}{2M-i}$$
and 
$$\frac{{2M-v\choose M}}{{2M\choose M}}=\frac{M}{2M}\frac{M-1}{2M-1}\ldots\frac{M-v+1}{2M-v+1} = \prod_{i=0}^{v-1}\frac{M-i}{2M-i}.$$
Thus

\[Z_3(M+1,M+1)=2\cdot\left(\frac{M+1}{2M+1}\right)^{2}\sum_{u=1,v=1}^{M}\frac{1}{u+v}\prod_{i=0}^{u-1}\frac{M-i}{2M-i}\prod_{i=0}^{v-1}\frac{M-i}{2M-i}\]

It follows that $w_{M+1,M+1} \leq \frac{M+1}{2M+1}$ is equivalent to the condition

\[2\left(\frac{M+1}{2M+1}\right)^2 \! \! \sum_{s=M+1}^{2M}\frac{1}{s}+2\left(\frac{M+1}{2M+1}\right)^{2}\! \! \sum_{u=1,v=1}^{M}\frac{1}{u+v}\prod_{i=0}^{u-1}\frac{M-i}{2M-i}\prod_{i=0}^{v-1}\frac{M-i}{2M-i} \leq \frac{M+1}{2M+1}.\]
This is equaivalent to
\begin{equation} 
\sum_{s=M+1}^{2M}\frac{1}{s}+\sum_{u=1,v=1}^{M}\frac{1}{u+v}\prod_{i=0}^{u-1}\frac{M-i}{2M-i}\prod_{i=0}^{v-1} \frac{M-i}{2M-i} \leq \frac{2M+1}{2M+2}.\label{principal}
\end{equation}
We have
\[\sum_{u=1,v=1}^{M}\frac{1}{u+v}\prod_{i=0}^{u-1}\frac{M-i}{2M-i}\prod_{i=0}^{v-1} \frac{M-i}{2M-i} \]$$=\frac{1}{8}+2\sum_{v=2}^M\frac{1}{v+1}\cdot\frac{1}{2}\cdot\frac{1}{2}\cdot\prod_{i=1}^{v-1}\cdot\frac{M-i}{2M-i}+\frac{1}{4}\sum_{u,v=2}^M\frac{1}{u+v}\prod_{i=1}^{u-1}\frac{M-i}{2M-i}\prod_{i=1}^{v-1}\frac{M-i}{2M-i}$$
$$=\frac{1}{8}+\frac{1}{2}\left(\frac{1}{3}\cdot\frac{M-1}{2M-1}+\frac{M-1}{2M-1}\sum_{v=3}^M\frac{1}{v+1}\prod_{i=2}^{v-1}\frac{M-i}{2M-i}\right)+$$
\[+\frac{1}{4}\! \left[\! \frac{1}{4}\! \left(\frac{M-1}{2M\! -\! 1}\right)^2 \! \! \! +\! 2\! \sum_{v=3}^M\frac{1}{v\! +\! 2}\cdot\frac{M-1}{2M\! -\! 1}\cdot\prod_{i=1}^{v-1}\frac{M-i}{2M\! -\! i}\! +\! \! \sum_{u,v=3}^M\frac{1}{u\! +\! v}\prod_{i=1}^{u-1}\frac{M-i}{2M\! -\! i}\prod_{i=1}^{v-1}\frac{M-i}{2M\! -\! i}\right]\]
$$=\frac{1}{8}+\frac{1}{2}\cdot\frac{M-1}{2M-1}\left(\frac{1}{3}+\sum_{v=3}^M\frac{1}{v+1}\prod_{i=2}^{v-1}\frac{M-i}{2M-i}\right)+\frac{1}{16} \left( \frac{M-1}{2M-1}\right) ^2+$$
\[+\frac{1}{4}\! \left[\! 2\! \left(\frac{M-1}{2M\! -\! 1}\right)^2\!  \sum_{v=3}^M\frac{1}{v\! +\! 2}\cdot\prod_{i=2}^{v-1}\! \frac{M-i}{2M\! -\! i}\! +
\! \left(\frac{M-1}{2M\! -\! 1}\right)^2\! \! \sum_{u,v=3}^M\frac{1}{u\! +\! v}\prod_{i=2}^{u-1}\frac{M-i}{2M\! -\! i}\prod_{i=1}^{v-1}\frac{M-i}{2M\! -\! i}\right]\]
$$\leq \frac{1}{8}+\frac{1}{2}\cdot\frac{M-1}{2M-1}\cdot\left(\frac{1}{3}+\sum_{v=3}^\infty\frac{1}{v+1}\left(\frac{M-2}{2M-2}\right)^{v-2}\right)+$$
$$+\frac{1}{4}\cdot\left(\frac{M-1}{2M-1}\right)^2\cdot\left[\frac{1}{4}+2\sum_{v=3}^\infty\frac{1}{v+2}\left(\frac{M-2}{2M-2}\right)^{v-2}+\sum_{u,v=3}^\infty\frac{1}{u+v}\left(\frac{M-2}{2M-2}\right)^{u+v-4}\right].$$
Setting $x=\frac{M-2}{2M-2}$, this expression is given by 
$$\frac{1}{8}+\frac{1}{2}\cdot\frac{M-1}{2M-1}\left(\frac{1}{3}+\frac{1}{x^3}\sum_{v=3}^\infty\frac{1}{v+1}x^{v+1}\right)+$$
$$+\frac{1}{4}\cdot\left(\frac{M-1}{2M-1}\right)^2\cdot\left[\frac{1}{4}+2\frac{1}{x^4}\sum_{v=3}^\infty\frac{1}{v+2}x^{v+2}+\frac{1}{x^4}\sum_{u,v=3}^\infty\frac{1}{u+v}x
^{u+v}\right]$$
$$=\frac{1}{8}+\frac{1}{2}\cdot\frac{M-1}{2M-1}\left(\frac{1}{3}+\frac{1}{x^3}\sum_{k=4}^\infty\frac{1}{k}x^{k}\right)+$$
$$+\frac{1}{4}\cdot\left(\frac{M-1}{2M-1}\right)^2\cdot\left[\frac{1}{4}+\frac{2}{x^4}\sum_{k=5}^\infty\frac{1}{k}x^{k}+\frac{1}{x^4}\sum_{k=6}^\infty\frac{k-5}{k}x^{k}\right]$$
$$=\frac{1}{8}+\frac{1}{2}\cdot\frac{M-1}{2M-1}\left(\frac{1}{3}+\frac{1}{x^3}\left(\sum_{k=1}^\infty\frac{x^k}{k}-\frac{x^3}{3}-\frac{x^2}{2}-x\right)\right)+$$
$$+\frac{1}{4}\cdot\left(\frac{M-1}{2M-1}\right)^2\cdot\Bigg[\frac{1}{4}+\frac{2}{x^4}\left(\sum_{k=1}^\infty\frac{x^k}{k}-\frac{x^4}{4}-\frac{x^3}{3}-\frac{x^2}{2}-x\right)+$$
$$+\frac{1}{x^4}\left(\sum_{k=6}^\infty x^k-5\left(\sum_{k=1}^\infty\frac{x^k}{k}-\frac{x^5}{5}-\frac{x^4}{4}-\frac{x^3}{3}-\frac{x^2}{2}-x\right)\right)\Bigg]$$
$$=\frac{1}{8}+\frac{1}{2}\cdot\frac{M-1}{2M-1}\left(\frac{1}{3}+\frac{1}{x^3}\left(-\ln(1-x)-\frac{x^3}{3}-\frac{x^2}{2}-x\right)\right)+$$
$$+\frac{1}{4}\cdot\left(\frac{M-1}{2M-1}\right)^2\cdot\Bigg[\frac{1}{4}+\frac{2}{x^4}\left(-\ln(1-x)-\frac{x^4}{4}-\frac{x^3}{3}-\frac{x^2}{2}-x\right)+$$
$$+\frac{1}{x^4}\left(\frac{x^6}{1-x}-5\left(-\ln(1-x)-\frac{x^5}{5}-\frac{x^4}{4}-\frac{x^3}{3}-\frac{x^2}{2}-x\right)\right)\Bigg]$$
$$=\frac{1}{8}-\frac{1}{2}\cdot\frac{M-1}{2M-1}\cdot\frac{1}{x}\left(\frac{1}{2}+\frac{1}{x}\right)+$$
$$+\frac{1}{4}\left(\frac{M-1}{2M-1}\right)^2\left(1+x+\frac{1}{x}+\frac{3}{2}\cdot\frac{1}{x^2}+\frac{3}{x^3}+\frac{x^2}{1-x}\right)+$$
$$-\ln(1-x)\cdot\frac{1}{2}\cdot\frac{M-1}{2M-1}\cdot\frac{1}{x^4}\left(x-\frac{3}{2}\cdot\frac{M-1}{2M-1}\right)$$
$$=\frac{1}{8}-\frac{1}{2}\cdot\frac{M-1}{2M-1}\cdot\frac{1}{x}\left(\frac{1}{2}+\frac{1}{x}\right)+$$
$$+\frac{1}{4}\left(\frac{M-1}{2M-1}\right)^2\left(1+x+\frac{1}{x}+\frac{3}{2}\cdot\frac{1}{x^2}+\frac{3}{x^3}+\frac{x^2}{1-x}\right)+$$
$$+\left(\ln(2)-\ln\left(1+\frac{1}{M-1}\right)\right)\cdot\frac{1}{2}\cdot\frac{M-1}{2M-1}\cdot\frac{1}{x^4}\left(x-\frac{3}{2}\cdot\frac{M-1}{2M-1}\right)$$
$$=\frac{1}{8}-\frac{1}{2}\cdot\frac{M-1}{2M-1}\cdot\frac{1}{x}\left(\frac{1}{2}+\frac{1}{x}\right)+$$
$$+\frac{1}{4}\left(\frac{M-1}{2M-1}\right)^2\left(1+x+\frac{1}{x}+\frac{3}{2}\cdot\frac{1}{x^2}+\frac{3}{x^3}+\frac{x^2}{1-x}\right)+$$
$$+\left(-\ln(2)+\ln\left(1+\frac{1}{M-1}\right)\right)\cdot\frac{1}{4}\cdot\frac{M^2-M+1}{x^4(2M-1)^2}.$$
Using the Taylor expansion of the logarithm expansion up to the third term, it follows that this expression is bounded above by
$$\frac{1}{8}-\frac{1}{2}\cdot\frac{M-1}{2M-1}\cdot\frac{1}{x}\left(\frac{1}{2}+\frac{1}{x}\right)+$$
$$+\frac{1}{4}\left(\frac{M-1}{2M-1}\right)^2\left(1+x+\frac{1}{x}+\frac{3}{2}\cdot\frac{1}{x^2}+\frac{3}{x^3}+\frac{x^2}{1-x}\right)+$$
$$+\left(-\ln(2)+\frac{1}{M-1}-\frac{1}{2(M-1)^2}+\frac{1}{3(M-1)^3}\right)\cdot\frac{1}{4}\cdot\frac{M^2-M+1}{x^4(2M-1)^2}.$$
Hence, $Z_3(M+1,M+1)$ has been bounded from above by a rational function. To obtain an upper bound on $2Z_1(M+1,M+1)$, note that
$$\left(\frac{M+1}{2M+1}\right)^{-2} Z_1(M+1,M+1)=\sum_{s=M+1}^{2M}\frac{1}{s}=
\sum_{s=M+1}^{2M}\frac{1}{s}-\ln(2)+\ln(2)$$
$$=\sum_{s=M+1}^{2M}\frac{1}{s}-\int_{M}^{2M}\frac{1}{t}\,dt+\ln(2)=
\ln(2)-\sum_{k=M}^{2M-1}\left(\int_{k}^{k+1}\frac{1}{t}\,dt-\frac{1}{k+1}\right)$$
$$=\ln(2)-\sum_{k=M}^{2M-1}\left(\ln\left(\frac{k+1}{k}\right)-\frac{1}{k+1}\right)=
\ln(2)-\sum_{k=M}^{2M-1}\left(\ln\left(1+\frac{1}{k}\right)-\frac{1}{k+1}\right)$$
$$=\ln(2)-\sum_{k=M}^{2M-1}\left(\left(\frac{1}{k}-\frac{1}{2k^2}+\frac{1}{3k^3}-\ldots\right)-\frac{1}{k+1}\right)$$
$$\leq \ln(2)-\sum_{k=M}^{2M-1}\left(\left(\frac{1}{k}-\frac{1}{2k^2}\right)-\frac{1}{k+1}\right)=
\ln(2)-\left(\frac{1}{M}-\frac{1}{2M}\right)+\frac{1}{2}\sum_{k=M}^{2M-1}\frac{1}{k^2}$$
$$\leq \ln(2)-\frac{1}{2M}+\frac{1}{2}\int_{M-1}^{2M-1}\frac{1}{t^2}\,dt=
\ln(2)-\frac{1}{2M}+\frac{1}{2}\left(\frac{1}{M-1}-\frac{1}{2M-1}\right).$$

Thus inequality (\ref{principal}) is satisfied when the following inequality holds
$$\frac{1}{8}-\frac{1}{2}\cdot\frac{M-1}{2M-1}\cdot\frac{1}{x}\left(\frac{1}{2}+\frac{1}{x}\right)+$$
$$+\frac{1}{4}\left(\frac{M-1}{2M-1}\right)^2\left(1+x+\frac{1}{x}+\frac{3}{2}\cdot\frac{1}{x^2}+\frac{3}{x^3}+\frac{x^2}{1-x}\right)+$$
$$+\left(-\ln(2)+\frac{1}{M-1}-\frac{1}{2(M-1)^2}+\frac{1}{3(M-1)^3}\right)\cdot\frac{1}{4}\cdot\frac{M^2-M+1}{x^4(2M-1)^2}+$$
$$+\ln(2)-\frac{1}{2M}+\frac{1}{2}\left(\frac{1}{M-1}-\frac{1}{2M-1}\right)\leq\frac{2M+1}{2M+2}.$$
Rearraging this inequality leads to

$$(384\ln(2)-264)M^8+(-2040\ln(2)+1401)M^7+(3840\ln(2)-2744)M^6+$$
$$+(-2088\ln(2)+2347)M^5+(-2400\ln(2)-1544)M^4+$$
$$+(3840\ln(2)+2956)M^3+(-1824\ln(2)-3824)M^2+(288\ln(2)+2032)M-384\geq0.$$

It is a matter of simple checking that 
$$(384\ln(2)-264)M^8+(-2040\ln(2)+1401)M^7+(3840\ln(2)-2744)M^6\geq0$$
for all $M\geq10$,
and that
$$(-2088\ln(2)+2347)M^5+(-2400\ln(2)-1544)M^4\geq0$$
and
$$(3840\ln(2)+2956)M^3+(-1824\ln(2)-3824)M^2\geq0$$
and 
$$(288\ln(2)+2032)M-384\geq0$$
for all nonnegative $M$.

Inequality (\ref{principal}) can be checked ''by hand''  for $M\leq 10$ . 

\newpage

\begin{table}[h]
\caption{Description of the notation for the most often used variables, functions and sets.} \label{notation}
\begin{tabular*}{0.9\linewidth}{@{\extracolsep{\fill}} ll @{}} \toprule
Notation & Description \\ \midrule
$N$ & Total number of objects of a given type \\ 
$i, j$ & Number of objects DM has rejected of type 1 and of type 2, respectively \\ 
$n=i+j$ & Number of objects rejected so far \\
$\tau_m$ & Total number of objects that have appeared when the $m$-th object  \\ 
& of
relative rank 1 appears (of either type) \\
$\alpha_{i,j}$ & Probability that the next object will be of type 1 \\ 
$v_{i,j}$ & Optimal expected reward of DM given the number of objects of each type rejected \\ 
$i^{\ast}_N (j)$ & Given that $j$ objects of type 2 have appeared, the DM should accept a candidate \\ & of type 1 when at least $i^{\ast}_N (j)$ objects
of type 1 have appeared. \\ 
$w_{i,j}$ & Expected reward of DM given that the next candidate is accepted. \\
$i^{\dagger}_N$ & $\min\{ i: v_{i,i}\leq \frac{i}{N}\}$ \\
$i^{\bullet}_N$ & $\min\{ i: w_{i,i} \leq \frac{i}{N}\}$ \\ 
$j^{\circ} (i)$ & Given that $i$ objects of type 1 have appeared, the next candidate to 
 appear should \\ & always be accepted if there have been at least $j^{\circ}(i)$ objects of type 2 \\ 
$C_0^N$ & Set of $(i,j)$ for which an object should never be accepted \\ 
$C_1^N$ & Set of $(i,j)$ for which only candidates of type 1 should be accepted \\ 
$C_2^N$ & Set of $(i,j)$ for which only candidates of type 2 should be accepted \\
$C_3^N$ & Set of $(i,j)$ for which candidates of either type should be accepted \\
$B_0^N$ & Set of $(i,j)$ for which the next object to appear should never be accepted \\ 
$B_1^N$ & Set of $(i,j)$ for which the next object to appear should only be accepted \\
& if it is a candidate of type 1 \\ 
$B_2^N$ & Set of $(i,j)$ for which the next object to appear should only be accepted \\
& if it is a candidate of type 2 \\ 
$B_3^N$ & Set of $(i,j)$ for which the next object to appear should be accepted if it is a candidate \\
\bottomrule
\end{tabular*}
\end{table}

\end{document}